\documentclass[reqno]{amsart}
\usepackage{amssymb}
\usepackage{graphicx}

\usepackage[usenames, dvipsnames]{color}
\usepackage{hyperref}
\usepackage{verbatim}
\usepackage{mathrsfs}
\usepackage{bm}
\usepackage{cite}

\numberwithin{equation}{section}

\newtheorem{theorem}{Theorem}[section]

\theoremstyle{definition}
\newtheorem{remark}[theorem]{Remark}

\theoremstyle{definition}

\theoremstyle{definition}

\makeatletter
\def\dashint{\operatorname%
{\,\,\text{\bf-}\kern-.98em\DOTSI\intop\ilimits@\!\!}}
\makeatother

\def\\det{\text{det}}

\def\.5{\frac{1}{2}}

\newcommand{\RN}[1]{%
  \textup{\uppercase\expandafter{\romannumeral#1}}%
}

\newcommand{\dist}{\text{dist}}

\renewcommand{\epsilon}{\varepsilon}

\newcounter{marnote}

\begin{document}
\title[Blow-up analysis of hydrodynamic forces ]{Blow-up analysis of hydrodynamic forces exerted on two adjacent $M$-convex particles}

\author[H.G. Li]{Haigang Li}
\address[H.G. Li]{School of Mathematical Sciences, Beijing Normal University, Beijing 100875, China. }
\email{hgli@bnu.edu.cn}

\author[X.T. Wang]{Xueting Wang}

\address[X.T. Wang]{School of Mathematical Sciences, Beijing Normal University, Beijing 100875, China. }
\email{201721130038@mail.bnu.edu.cn}

\author[Z.W. Zhao]{Zhiwen Zhao}

\address[Z.W. Zhao]{1. School of Mathematical Sciences, Beijing Normal University, Beijing 100875, China. }
\address{2. Bernoulli Institute for Mathematics, Computer Science and Artificial Intelligence, University of Groningen, PO Box 407, 9700 AK Groningen, The Netherlands.}
\email{zwzhao@mail.bnu.edu.cn, Corresponding author.}


\date{\today} 




\begin{abstract}
In a viscous incompressible fluid, the hydrodynamic forces acting on two close-to-touch rigid particles in relative motion always become arbitrarily large, as the interparticle distance parameter $\varepsilon$ goes to zero. In this paper we obtain asymptotic formulas of the hydrodynamic forces and torque in $2\mathrm{D}$ model and establish the optimal upper and lower bound estimates in $3\mathrm{D}$, which sharply characterizes the singular behavior of hydrodynamic forces. These results reveal the effect of the relative convexity between particles, denoted by index $m$, on the blow-up rates of hydrodynamic forces. Further, when $m$ degenerates to infinity, we consider the particles with partially flat boundary and capture that the largest blow-up rate of the hydrodynamic forces is $\varepsilon^{-3}$ both in 2D and 3D. We also clarify the singularities arising from linear motion and rotational motion, and find that the largest blow-up rate induced by rotation appears in all directions of the forces.
\end{abstract}

\maketitle



\section{Introduction}

For suspension problem of rigid particles in a viscous incompressible fluid, there are wide applications in nature and engineering, such as in geology, biology, fluid mechanics, chemical engineering and composite material manufacturing, see e.g. \cite{G1994,S1989,S2011,PT2009}. In this paper we consider a flow with low Reynolds number so that particle inertia, fluid inertia and Brownian motion may be neglected. When two rigid particles in a viscous incompressible fluid are located closely to each other, the hydrodynamic forces acting on particles in relative motion will appear blow-up. To better understand this high concentration, it is significantly important to give a quantitative characterization of the singularities of the hydrodynamic forces in terms of the distance, say $\varepsilon$, between particles.

In the past few decades, there has been a long list of literature on the study of the singularities of hydrodynamic forces. For small $\varepsilon$, Goldman, Cox and Brenner \cite{GCB1967I} developed a lubrication theory and obtained an asymptotic expansion, with blow-up rate $|\ln\varepsilon|$, for the force and torque exerted by the fluid on the sphere moving parallel to a solid plane wall. Cox and Brenner \cite{CB1967II} adopted the singular perturbation techniques to calculate the force experienced by a sphere moving perpendicularly to a solid plane wall and obtained the blow-up rate is $\varepsilon^{-1}$. Cox \cite{C1974} extended to the case of $2$-convex surfaces (corresponding to the case that $m=2$ in definition \eqref{ZAD01} below) and obtained the first leading singular term with rate $\varepsilon^{-1}$ and the second term with rate $|\ln\varepsilon|$ for the force and the torque. For more related works, see \cite{G1980,J1982,JO1984,NK1984,CB1989} and references therein. In these above studies, the blow-up rate $\varepsilon^{-1}$ has been accurately derived, while the terms with $|\ln\varepsilon|$ not yet fully understood. This disagreement is mainly due to the difference of approximation method and mechanisms to be dealt with. Some people consider translational motion rather than rotation, while others adopt an approximation for the surface of particles, which only contributes to capturing accurately the major blow-up rate $\varepsilon^{-1}$. In order to overcome this shortcoming caused by the asymptotic results for the rate $|\ln\varepsilon|$, Gorb \cite{G2016} developed a concise analytical method to derive all the blow-up rates of hydrodynamic forces for spherical particles in $3\mathrm{D}$.

Besides of the hydrodynamic forces, the effective viscosity is another quantity of interest, which describes the effective rheological properties of suspensions. Frankel and Acrivos \cite{FA1967} performed a formal asymptotic analysis of the effective viscosity in the narrow region between two neighbouring spherical inclusions. They considered the translation motions of particles and obtained the asymptotics in form of $C\varepsilon^{-1}+O(|\ln\varepsilon|)$. Subsequently, a periodic array of particles in a Newtonian fluid was investigated in \cite{NK1984}. Sierou and Brady \cite{SB2001,SB2002} investigated numerically the concentrated random suspensions by using accelerated Stokesian dynamics. It was observed that in some cases the effective viscosity is of order $\varepsilon^{-1}$, while in other cases it is of order $|\ln\varepsilon|$. This means that for generic random suspensions, the asymptotic behavior of the effective viscosity defined by the global dissipation rate cannot be determined by the local dissipation rate in a single gap. For more related numerical results, we refer to papers \cite{M1999,GPHJP2001,M1997,DM2003,SC2012}. Given certain boundary conditions and an array of particles in 3D, it was shown that the blow-up rate $|\ln\varepsilon|$ becomes dominant with the degeneration of the leading term $O(\varepsilon^{-1})$. In fact, due to reduced analytical and computational complexity, 2D models are often used to describe 3D suspensions. Berlyand, Gorb and Novikov \cite{BGN2009} proved the validity of 2D models for 3D problems.

We would like to point out that Stokes problem is actually a limit of Lam\'e problem in linear elasticity when the Lam\'e coefficient $\lambda$ goes to infinity \cite{MS2012,AKKL2013}. This is also one motivation of this paper, because these high concentration phenomenon are quite similar. When the distance $\varepsilon$ between two stiff inclusions tends to zero, it is proved that the stress field blows up at rate $\varepsilon^{-1/2}$ in $2\mathrm{D}$ \cite{BLL2015,KY2019, Li2018}, and $(\varepsilon|\ln\varepsilon|)^{-1}$ in $3\mathrm{D}$ \cite{BLL2017,Li2018} for 2-convex inclusions. For general $m$-convex inclusions, it is also studied in \cite{HL2018,LZ2020}. The singularities of the stress field will disappear when $m$ goes to infinity. This is contrary to the concentration phenomena of the hydrodynamic forces occurring in suspensions of rigid particles as two stiff particles are close to touching. The analogous electric field concentration was also studied, see \cite{BLY2009,BT2013,AKL2005,Y2007,KLY2013,KLY2014,LLY2019} and the reference therein. Finally, it is worth mentioning that Ammari, Kang, Kim and Yu \cite{AKKY2020} recently showed that the blow-up rate of the fluid stress is $\varepsilon^{-\frac{1}{2}}$ in $2\mathrm{D}$.

In this paper, we will study the $m$-convex particles ($m\geq2$, the definition given by \eqref{ZAD01} below) both in $2\mathrm{D}$ and $3\mathrm{D}$, and establish the asymptotic formula of hydrodynamic forces to reveal their explicit dependence on the convexity of particles. It shows that the singularity of these forces strengthens as $m$ increases. When $m$ further degenerates to $\infty$, the shapes of particles becomes flat. In this case, we capture the largest blow-up rate $\varepsilon^{-3}$. Additionally, our results classify the singularities according to linear motion and rotational motion, respectively. We show that the largest blow-up rate generated by rotation appears in all directions of the forces.

The rest of this paper is organized as follows. In Section 2, we describe the problem and state our main results, including Theorems \ref{thm1} and \ref{thm3} for $3\mathrm{D}$ and Theorem \ref{thm2} for $2\mathrm{D}$. In Section 3 we present the approximation method in $3\mathrm{D}$. The proofs of Theorems \ref{thm1} and \ref{thm3} are given in Sections 4 and 5, respectively. Finally, the approximation method under the case of $2\mathrm{D}$ and the proof of Theorem \ref{thm2} are shown in Sections 6 and 7, respectively.

\section{Formulation of the problem and main results}
\subsection{Problem formulation}

Let $\Omega\subseteq\mathbb{R}^{3}$ be a smooth region occupied by an incompressible fluid with viscosity $\mu$, which contains two identical convex rigid particles $D_{1}$ and $D_{2}$ with $C^{2}$ boundaries. We assume that $D_{1}$ and $D_{2}$ are symmetric with respect to $x_{3}$-axis and the plane $\{x_{3}=0\}$, and their centroids both lie on the $x_{3}$-axis. Denote $x'=(x_{1},x_{2})$. Suppose the interparticle distance between them is $\varepsilon$, a sufficiently small positive constant and achieves at two points, $(0',\frac{\varepsilon}{2})\in\partial D_{1}$ and $(0',-\frac{\varepsilon}{2})\in\partial D_{2}$. Assume that the centers of mass of particles $D_{1}$ and $D_{2}$ are at $\left(0',\pm\frac{\varepsilon}{2}\pm R\right)$, respectively, for a fixed constant $R>0$. Fix a small constant $0<r<R$ such that the portions of $\partial D_{1}$ and $\partial D_{2}$ near the origin can be expressed, respectively, as
\begin{align}\label{ZAD01}
x_{3}=\frac{\varepsilon}{2}+\frac{a}{2}|x'|^{m}\quad\mathrm{and}\quad x_{3}=-\frac{\varepsilon}{2}-\frac{b}{2}|x'|^{m},\quad m\geq2,\;\quad \mathrm{for}\;|x'|<r.
\end{align}
We call these particles {\em $m$-convex particles}. For simplicity, we assume that $a=b=1$ in the following. When particle $D_{2}$ remains stationary, we let particle $D_{1}$, respectively, move at constant linear velocity, $\mathbf{U}$, and constant angular velocity, $\boldsymbol{\omega}$, with
\begin{align*}
 \mathbf{U}=
U_{1}\mathbf{e}_{1}+U_{2}\mathbf{e}_{2}+U_{3}\mathbf{e}_{3},
\end{align*}
and
\[
\ \boldsymbol{\omega}=\omega_{1}\mathbf{e}_{1}+\omega_{2}\mathbf{e}_{2}+\omega_{3}\mathbf{e}_{3},
\]
where $\{\mathbf{e}_{1},\mathbf{e}_{2},\mathbf{e}_{3}\}$ represents the standard basis of three-dimensional Euclidean space. Here we consider a stationary flow regime and assume that $U_{i}$, $\omega_{i}$, $\dist(D_{1}\cup D_{2},\partial\Omega)$, and the diameters of $D_{1}$ and $D_{2}$ are constants independent of $\varepsilon$.

Let $\mathbf{u}$ be the velocity of the fluid and $p$ be the pressure. Denote the rate of strain tensor by
\begin{equation*}
\mathbf{D(u)}=\frac{1}{2}(\nabla\mathbf{u}+\nabla\mathbf{u}^{T}),
\end{equation*}
and the fluid stress by
\begin{equation*}
 \boldsymbol{\sigma}\mathbf{(u)}=2\mu\mathbf{D(u)}-p\mathbf{I}.
\end{equation*}
Set $\mathbf{x}_{D_{1}}$ to be the center of mass of $D_{1}$. We consider the following boundary value problem:
\begin{align}\label{11.3}
\begin{cases}
(\mathrm{a})~\nabla\cdot\boldsymbol{\sigma}\mathbf{(u)}=\mathbf{0}&~~~~~\mathrm{in}~\Omega,\\
(\mathrm{b})~~~~~~~\nabla\cdot\mathbf{u}=0&~~~~~\mathrm{in}~\Omega,\\
(\mathrm{c})~~~~~~~~~~~~~~\mathbf{u}=\mathbf{U}+\boldsymbol{\omega}\times(\mathbf{x}-\mathbf{x}_{D_{1}})&~~~~~\mathrm{on}~\partial D_{1},\\
(\mathrm{d})~~~~~~~~~~~~~~\mathbf{u}=\mathbf{0}&~~~~~\mathrm{on}~\partial D_{2},\\
(\mathrm{e})~\boldsymbol{\sigma}\mathbf{(u)n}_{\partial\Omega}=\mathbf{0}&~~~~~\mathrm{on}~\partial\Omega,
\end{cases}
\end{align}
where $\mathbf{u}\in \mathbf{H}^{1}(\Omega)$ and $\mathbf{n}_{\partial\Omega}$ denotes the unit outer normal to $\partial\Omega$. $(\ref{11.3})(\mathrm{b})$ provides the incompressibility for the stationary Stokes flow. $(\ref{11.3})(\mathrm{c})$ gives the no-slip condition of the fluid velocity on the surface of $D_{1}$. $(\ref{11.3})(\mathrm{d})$ keeps particle $D_{2}$ motionless. $(\ref{11.3})(e)$ ensures that the fluid does not flow out of the external boundary $\partial\Omega$.

Throughout this paper, for simplicity, we denote
\begin{equation}\label{PLG01}
\mathbf{x}-\mathbf{x}_{D_{1}}:=\boldsymbol{\nu}=(\nu_{1},...,\nu_{n}),\quad n=2,3.
\end{equation}
Denote by
\begin{equation}\label{11.4}
\mathbf{F}=(F_{1},F_{2},F_{3})=\int_{\partial D_{1}}\boldsymbol{\sigma}\mathbf{(u)n}\,dS
\end{equation}
{\em the hydrodynamic force} acting on $D_{1}$ and
\begin{equation}\label{11.5}
\mathbf{T}=(T_{1},T_{2},T_{3})=\int_{\partial D_{1}}\boldsymbol{\nu}\times\boldsymbol{\sigma}\mathbf{(u)n}\,dS
\end{equation}
{\em the hydrodynamic torque}, where $\mathbf{n}$ is the unit outer normal to $\partial D_{1}$ and $\boldsymbol{\nu}$ is defined in \eqref{PLG01}. We use the notation
\begin{align*}
\Gamma_{ij}^{(m)}&:=
\begin{cases}
\frac{1}{m}\Gamma\left(i-\frac{j}{m}\right)\Gamma\left(\frac{j}{m}\right),&i\neq\frac{j}{m},\\
\frac{1}{m},&i=\frac{j}{m},
\end{cases}\quad\mathrm{for}\;ij\in\{11,12,13,23,24,33,34,35,36\},
\end{align*}
where $\Gamma(s)=\int^{+\infty}_{0}x^{s-1}e^{-x}dx$, $s>0$ is the Gamma function. Without loss of generality, we assume
$$U_{3}\leq0,~\omega_{i}\geq0,\quad i=1,2,3.$$
\subsection{Main results}

Our main results are as follows. We separate it into two parts: $m=2$ and $m>2$.

\begin{theorem}\label{thm1}
Assume that $D_{1},D_{2}\subset \Omega\subseteq\mathbb{R}^{3}$ are defined as above and \eqref{ZAD01} holds. Let $\mathbf{u}\in \mathbf{H}^{1}(\Omega)$ be the solution of problem (\ref{11.3}). Then, for a sufficiently small $\varepsilon,$ we have the following asymptotic expansions:

$(i)$ for $m=2$,
\begin{align*}
\frac{\omega_{2}r^{2}\alpha^{(2)}_{34}}{4\varepsilon}\leq& F_{1}+(U_{1}-\omega_{2}R)\alpha_{12}^{(2)}|\ln\varepsilon|\leq \frac{2\omega_{2}r^{2}\alpha_{34}^{(2)}}{\varepsilon},\\
\frac{\omega_{1}r^{2}\alpha_{34}^{(2)}}{4\varepsilon}\leq -&F_{2}-(U_{2}+\omega_{1}R)\alpha_{12}^{(2)}|\ln\varepsilon|\leq\frac{2\omega_{1}r^{2}\alpha_{34}^{(2)}}{\varepsilon},\\
\frac{(\omega_{1}+\omega_{2})r\alpha_{34}^{(2)}}{\varepsilon}\leq& F_{3}+2U_{3}\alpha_{34}^{(2)}\frac{1}{\varepsilon}\leq\frac{2(\omega_{1}+\omega_{2})r\alpha_{34}^{(2)}}{\varepsilon},
\end{align*}
and
\begin{align*}
\frac{\omega_{1}r^{2}\beta_{1}^{(2)}}{\varepsilon}\leq\,&T_{1}+R(U_{2}+\omega_{1}R)\alpha_{12}^{(2)}|\ln\varepsilon|\leq\frac{\omega_{1}r^{2}\beta_{2}^{(2)}}{\varepsilon},\\
\frac{\omega_{2}r^{2}\beta_{1}^{(2)}}{\varepsilon}\leq\,&T_{2}-R(U_{1}-\omega_{2}R)\alpha_{12}^{(2)}|\ln\varepsilon|\leq\frac{\omega_{2}r^{2}\beta_{2}^{(2)}}{\varepsilon},\\
&T_{3}=0;
\end{align*}

$(ii)$ for $m>2$,
\begin{align*}
\frac{\omega_{2}2^{-m}r^{m}\alpha_{34}^{(m)}}{\varepsilon^{3-{4}/{m}}}\leq& F_{1}+\frac{(U_{1}-\omega_{2}R)\alpha_{12}^{(m)}}{\varepsilon^{1-{2}/{m}}}-\frac{\omega_{2}\alpha_{34}^{(m)}}{\varepsilon^{2-{4}/{m}}}\leq \frac{\omega_{2}2^{\frac{m}{2}}r^{m}\alpha_{34}^{(m)}}{\varepsilon^{3-{4}/{m}}},\\
\frac{\omega_{1}2^{-m}r^{m}\alpha_{34}^{(m)}}{\varepsilon^{3-{4}/{m}}}\leq -&F_{2}-\frac{(U_{2}+\omega_{1}R)\alpha_{12}^{(m)}}{\varepsilon^{1-{2}/{m}}}-\frac{\omega_{1}\alpha_{34}^{(m)}}{\varepsilon^{2-{4}/{m}}}\leq\frac{\omega_{1}2^{\frac{m}{2}}r^{m}\alpha_{34}^{(m)}}{\varepsilon^{3-{4}/{m}}},\\
\frac{(\omega_{1}+\omega_{2})r\alpha_{34}^{(m)}}{\varepsilon^{3-{4}/{m}}}\leq& F_{3}+\frac{2U_{3}\alpha_{34}^{(m)}}{\varepsilon^{3-{4}/{m}}}\leq\frac{2(\omega_{1}+\omega_{2})r\alpha_{34}^{(m)}}{\varepsilon^{3-{4}/{m}}},
\end{align*}
and
\begin{align*}
\frac{\omega_{1}r^{2}\beta_{1}^{(m)}}{\varepsilon^{3-{4}/{m}}}\leq\,& T_{1}+\frac{R(U_{2}+\omega_{1}R)\alpha_{12}^{(m)}}{\varepsilon^{1-{2}/{m}}}\leq\frac{\omega_{1}r^{2}\beta_{2}^{(m)}}{\varepsilon^{3-{4}/{m}}},\\
\frac{\omega_{2}r^{2}\beta_{1}^{(m)}}{\varepsilon^{3-{4}/{m}}}\leq\,& T_{2}-\frac{R(U_{1}-\omega_{2}R)\alpha_{12}^{(m)}}{\varepsilon^{1-{2}/{m}}}\leq\frac{\omega_{2}r^{2}\beta_{2}^{(m)}}{\varepsilon^{3-{4}/{m}}},\\
&T_{3}=0;
\end{align*}
where
\begin{align*}
\alpha_{12}^{(m)}=&2\pi\mu\Gamma_{12}^{(m)},\;\beta^{(m)}_{1}=\frac{3}{16}\pi\mu\Gamma_{34}^{(m)}(1-2^{\frac{m}{2}+2}Rr^{m-2}+2^{-2m}r^{2m-2}),\\ \alpha_{34}^{(m)}=&\frac{3}{2}\pi\mu\Gamma_{34}^{(m)},\;\beta^{(m)}_{2}=3\pi\mu\Gamma_{34}^{(m)}(1-2^{-m}Rr^{m-2}+2^{m-2}r^{2m-2}),
\end{align*}
and we omit the remainder terms of order $O(1)$, depending only on $\mu$, $R$, $r$, $U_{i}$ and $\omega_{i}$ but not on $\varepsilon$, in the estimates above of each component of $\mathbf{F}$ and $\mathbf{T}$, respectively.
\end{theorem}
\begin{remark}
From Theorem \ref{thm1}, we find that when $m>2$, there are three blow-up terms with rates $\frac{1}{\varepsilon^{1-{2}/{m}}}$, $\frac{1}{\varepsilon^{2-{4}/{m}}}$, and $\frac{1}{\varepsilon^{3-{4}/{m}}}$ in terms of the hydrodynamic force $\mathbf{F}$. These singularities of hydrodynamic forces all increase as $m$ increases, meaning that the relative convexity of surfaces weakens. Moreover, the largest blow-up rate $\frac{1}{\varepsilon^{3-{4}/{m}}}$ caused by rotation appears in all directions, while the rate $\frac{1}{\varepsilon^{3-{4}/{m}}}$ induced by the linear motion only appears in the direction of $x_{3}$-axis, which indicates that the singularities created by rotational motion hold dominant position in $x_{1}$-axis and $x_{2}$-axis rather than linear motion. Finally, we see from the coefficients of the blow-up rates of the hydrodynamic forces in Theorem \ref{thm1} that the singular effect of the geometry parameters $r$ and $R$ is embodied with rotational motion.
\end{remark}
\begin{remark}
If $\boldsymbol{\omega}=\mathbf{0}$ in Theorem \ref{thm1}, then

$(i)$ for $m=2$,
\begin{align*}
\mathbf{F}=&-\pi\mu
\left(
\begin{array}{c}
          |\ln\varepsilon|U_{1}\\
          |\ln\varepsilon|U_{2}\\
          \frac{3}{\varepsilon}U_{3}
\end{array}
 \right)+O(\mathbf{1}),\quad
\mathbf{T}=\pi\mu R|\ln\varepsilon|
\left(
\begin{array}{c}
          -U_{2}\\
          U_{1}\\
          0
\end{array}
 \right)+O(\mathbf{1});
\end{align*}

$(ii)$ for $m>2$,
\begin{align*}
\mathbf{F}=&-\frac{2\pi\mu}{\varepsilon^{1-{2}/{m}}}
\left(
\begin{array}{c}
          \Gamma_{12}^{(m)}U_{1}\\
          \Gamma_{12}^{(m)}U_{2}\\
          \frac{3}{2\varepsilon^{2-{2}/{m}}}\Gamma_{34}^{(m)}U_{3}
\end{array}
 \right)+O(\mathbf{1}),\quad
\mathbf{T}=\frac{2\pi\mu R\Gamma_{12}^{(m)}}{\varepsilon^{1-{2}/{m}}}
\left(
\begin{array}{c}
          -U_{2}\\
          U_{1}\\
          0
\end{array}
 \right)+O(\mathbf{1}),
\end{align*}
where the remainder terms of order $O(\mathbf{1})$ depend only on $\mu$, $R$, $r$, $U_{1}$, $U_{2}$, but not on $\varepsilon$. The result of $m=2$ is consistent with that in \cite{G2016}. We extend to the case when $m>2$ and study the effect of the geometry of the particles on the hydrodynamic forces.
\end{remark}

\begin{remark}
Generally, if
\begin{align*}
x_{3}=\frac{\varepsilon}{2}+\frac{k_{1}}{2}|x_{1}|^{m}+\frac{k_{2}}{2}|x_{2}|^{m},~\mathrm{and}~ x_{3}=-\frac{\varepsilon}{2}-\frac{k_{1}}{2}|x_{1}|^{m}-\frac{k_{2}}{2}|x_{2}|^{m},~\;\mathrm{for}\;|x'|<r,
\end{align*}
where $k_{1},\,k_{2}>0$ and $m\geq2$.  By using the proof of Theorem \ref{thm1} below with minor modification, we are able to obtain that the coefficients of the leading terms have an explicit dependence on Gaussian curvature $k_{1}k_{2}$ in the form of $(k_{1}k_{2})^{-\frac{1}{m}}$, which implies that the smaller the principal curvatures, the greater the singularities of hydrodynamic forces.
\end{remark}

From Theorem \ref{thm1}, we see that the singularities of hydrodynamic forces increases as $m\rightarrow +\infty$, the surfaces of particles in the narrow region partially tend to be flat. Next, we consider the particles with partially flat boundary. Take a small constant $0<s<r$, independent of $\varepsilon$. Assume that the corresponding partial boundaries of $\partial D_{1}$ and $\partial D_{2}$ are formulated, respectively, as
\begin{align}\label{2.001}
x_{3}&=
\begin{cases}
\frac{\varepsilon}{2},&|x'|\leq s,\\
\frac{\varepsilon}{2}+\frac{1}{2}(|x'|-s)^{2},&s<|x'|\leq r,
\end{cases}\;\,
x_{3}=\begin{cases}
-\frac{\varepsilon}{2},&|x'|\leq s,\\
-\frac{\varepsilon}{2}-\frac{1}{2}(|x'|-s)^{2},&s<|x'|\leq r.
\end{cases}
\end{align}
\begin{theorem}\label{thm3}
Assume that $D_{1},D_{2}\subset\Omega\subseteq\mathbb{R}^{3}$ are defined as above, condition (\ref{2.001}) holds for $0<s<(\sqrt{2}-1)r$. Let $\mathbf{u}\in \mathbf{H}^{1}(\Omega)$ be the solution of problem (\ref{11.3}). Then, for a sufficiently small $\varepsilon$,
\begin{align*}
\omega_{2}\mathcal{B}_{1}(\varepsilon^{-3})\leq& F_{1}+\pi\mu(U_{1}-\omega_{2}R)\left(|\ln\varepsilon|+\frac{2s\Gamma^{(2)}_{11}}{\varepsilon^{1/2}}+\frac{s^{2}}{\varepsilon}\right)\leq \omega_{2}\mathcal{B}_{2}(\varepsilon^{-3}),\\
\omega_{1}\mathcal{B}_{1}(\varepsilon^{-3})\leq& -F_{2}-\pi\mu(U_{2}+\omega_{1}R)\left(|\ln\varepsilon|+\frac{2s\Gamma^{(2)}_{11}}{\varepsilon^{1/2}}+\frac{s^{2}}{\varepsilon}\right)\leq\omega_{1}\mathcal{B}_{2}(\varepsilon^{-3}),\\
\sum^{2}_{i=1}\omega_{i}\mathcal{C}_{1}(\varepsilon^{-3})\leq& F_{3}+3\pi\mu U_{3}\left(\frac{1}{2\varepsilon}+\frac{3s\Gamma_{21}^{(2)}}{\varepsilon^{3/2}}+\frac{4 s^{2}\Gamma_{32}^{(2)}}{\varepsilon^{2}}+\frac{2r^{2}s^{2}-s^{4}}{2\varepsilon^{3}}\right)\leq\sum^{2}_{i=1}\omega_{i}\mathcal{C}_{2}(\varepsilon^{-3}),
\end{align*}
and
\begin{align*}
\omega_{1}\mathcal{D}_{1}(\varepsilon^{-3})\leq\, &T_{1}+\pi\mu R(U_{2}+\omega_{1}R)\bigg(|\ln\varepsilon|+\frac{2s\Gamma^{(2)}_{11}}{\varepsilon^{1/2}}+\frac{s^{2}}{\varepsilon}\bigg)\leq\omega_{1}\mathcal{D}_{2}(\varepsilon^{-3}),\\
\omega_{2}\mathcal{D}_{1}(\varepsilon^{-3})\leq\,& T_{2}-\pi\mu R(U_{1}-\omega_{2}R)\bigg(|\ln\varepsilon|+\frac{2s\Gamma^{(2)}_{11}}{\varepsilon^{1/2}}+\frac{s^{2}}{\varepsilon}\bigg)\leq\omega_{2}\mathcal{D}_{2}(\varepsilon^{-3}),\\
& T_{3}=0,
\end{align*}
where
\begin{align*}
\mathcal{B}_{1}(\varepsilon^{-3})=&\frac{3}{16}\pi\mu(r-s)^{2}\left(\frac{1}{\varepsilon}+\frac{s^{4}}{\varepsilon^{3}}\right),\;\mathcal{B}_{2}(\varepsilon^{-3})=\frac{3}{4}\pi\mu(2r-s)^{2}\left(\frac{1}{\varepsilon}+\frac{2s^{4}}{\varepsilon^{3}}\right),\\
\mathcal{C}_{1}(\varepsilon^{-3})=&\frac{3}{4}\pi\mu(r+s)\left(\frac{1}{\varepsilon}+\frac{s^{4}}{\varepsilon^{3}}\right),\;\mathcal{C}_{2}(\varepsilon^{-3})=\frac{3}{2}\pi\mu r\left(\frac{1}{\varepsilon}+\frac{2s^{4}}{\varepsilon^{3}}\right),\\
\mathcal{D}_{1}(\varepsilon^{-3})=&\frac{3}{32}\pi\big[-4R(\sqrt{2}r-s)^{2}+(r+s)^{2}+2^{-4}(r-s)^{4}\big]\varepsilon^{-1}\\
&-\frac{3}{16}\pi s^{4}\big[4R(\sqrt{2}r-s)^{2}+(s-r)^{2}-2r^{2}-2^{-4}(r-s)^{4}\big]\varepsilon^{-3},\\
\mathcal{D}_{2}(\varepsilon^{-3})=&\frac{3}{8}\pi\big[-R(r-s)^{2}+4r^{2}+(\sqrt{2}r-s)^{4}\big]\varepsilon^{-1}\\
&-\frac{3}{16}\pi s^{4}\big[R(r-s)^{2}-4r^{2}+2s^{2}-(\sqrt{2}r-s)^{4}\big]\varepsilon^{-3},
\end{align*}
and we omit the remainder terms of order $O(1)$ depending only on $\mu$, $R$, $r$, $s$, $U_{i}$ and $\omega_{i}$ but not on $\varepsilon$ in the estimates above of each component of $\mathbf{F}$ and $\mathbf{T}$, respectively.

\end{theorem}
\begin{remark}
When $s$ tends to zero, we get the estimates in Theorem \ref{thm1} for $m=2$. Furthermore, we prove that the largest blow-up rate is $\varepsilon^{-3}$ when the flatness appears.
\end{remark}
\begin{figure}[htb]
\center{\includegraphics[width=0.45\textwidth]{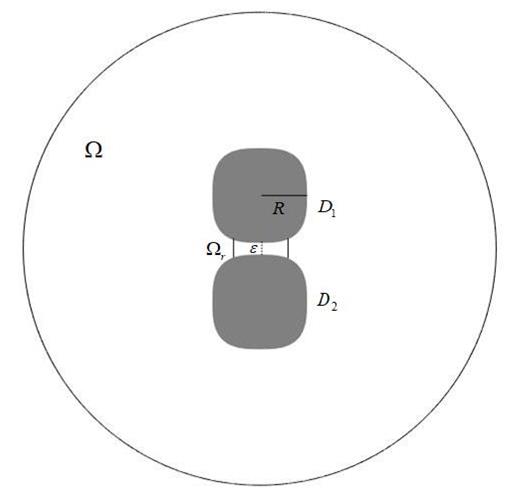}}
\caption{our domain $\Omega$}
\end{figure}

The corresponding $2\mathrm{D}$ problem is similar. Here we point out some differences between them. First, in $2\mathrm{D}$ the angular velocity and the hydrodynamic torque are scalars denoted by $\omega_{0}$ and $\mathrm{T}$, respectively. Second, assume that the centroid of particles $D_{1}$ and $D_{2}$ lie at $\left(0,\pm\frac{\varepsilon}{2}\pm R\right)$, respectively (see Figure 1), and $m>1$ in $2\mathrm{D}$. Finally, in contrast to the assumptions in $3\mathrm{D}$, in $2\mathrm{D}$ we need to replace condition $(\ref{11.3})(c)$ by
\begin{align*}
\mathbf{u}=\mathbf{U}+\omega_{0}\times\boldsymbol{\nu}=\mathbf{U}+\omega_{0}(\nu_{1}\mathbf{e}_{2}-\nu_{2}\mathbf{e}_{1}),\quad\mathrm{on}\;\partial D_{1},
\end{align*}
where $\omega_{0}$, a scalar, is angular velocity of particle $D_{1}$. That is, a two-dimensional ($2\mathrm{D}$) model is as follows:
\begin{align}\label{ZAD123}
\begin{cases}
(\mathrm{a})~\nabla\cdot\boldsymbol{\sigma}\mathbf{(u)}=\mathbf{0}&~~~~~\mathrm{in}~\Omega,\\
(\mathrm{b})~~~~~~~\nabla\cdot\mathbf{u}=0&~~~~~\mathrm{in}~\Omega,\\
(\mathrm{c})~~~~~~~~~~~~~~\mathbf{u}=\mathbf{U}+\omega_{0}\times\boldsymbol{\nu}&~~~~~\mathrm{on}~\partial D_{1},\\
(\mathrm{d})~~~~~~~~~~~~~~\mathbf{u}=\mathbf{0}&~~~~~\mathrm{on}~\partial D_{2},\\
(\mathrm{e})~\boldsymbol{\sigma}\mathbf{(u)n}_{\partial\Omega}=\mathbf{0}&~~~~~\mathrm{on}~\partial\Omega,
\end{cases}
\end{align}
where $\mathbf{u}\in \mathbf{H}^{1}(\Omega)$ and $\mathbf{n}_{\partial\Omega}$ represents the unit outer normal of $\partial\Omega$.

\begin{theorem}\label{thm2}
Assume that $D_{1},D_{2}\subset\Omega\subseteq\mathbb{R}^{2}$ are defined as above, condition \eqref{ZAD01} holds. Let $\mathbf{u}\in \mathbf{H}^{1}(\Omega)$ be the solution of problem \eqref{ZAD123}. Then, for a sufficiently small $\varepsilon$,

$(i)$~
if $1<m\leq\frac{3}{2}$,
\begin{align*}
\mathbf{F}=&-
\left(
\begin{array}{c}
          (U_{1}+\omega_{0}R)\alpha_{11}^{(m)}\frac{1}{\varepsilon^{1-1/m}}+\omega_{0}r^{m}\alpha_{33}^{(m)}\frac{1}{\varepsilon^{3-3/m}}\\\\
          2(2U_{2}-\omega_{0}r)\alpha_{33}^{(m)}\frac{1}{\varepsilon^{3-3/m}}\\
\end{array}
 \right)+O(\mathbf{1});
\end{align*}
if $m>\frac{3}{2}$,
\begin{align*}
\mathbf{F}=&-
\left(
\begin{array}{c}
          (U_{1}+\omega_{0}R)\alpha_{11}^{(m)}\frac{1}{\varepsilon^{1-1/m}}+\omega_{0}\alpha_{33}^{(m)}\frac{1}{\varepsilon^{2-3/m}}+\omega_{0}r^{m}\alpha_{33}^{(m)}\frac{1}{\varepsilon^{3-3/m}}\\\\
          2(2U_{2}-\omega_{0}r)\alpha_{33}^{(m)}\frac{1}{\varepsilon^{3-3/m}}\\
\end{array}
 \right)+O(\mathbf{1}).
\end{align*}

$(ii)$ if $1<m\leq\frac{3}{2}$,
\begin{align*}
\mathrm{T}=&-R(U_{1}+\omega_{0}R)\alpha_{11}^{(m)}\frac{1}{\varepsilon^{1-1/m}}+\omega_{0}r^{2}\beta^{(m)}\frac{1}{\varepsilon^{3-3/m}}+O(1);
\end{align*}
if $\frac{3}{2}<m<\frac{5}{3}$,
\begin{align*}
\mathrm{T}=&-\frac{R(U_{1}+\omega_{0}R)\alpha_{11}^{(m)}}{\varepsilon^{1-1/m}}-\frac{R\omega_{0}\alpha_{33}^{(m)}}{\varepsilon^{2-3/m}}+\frac{\omega_{0}r^{2}\beta^{(m)}}{\varepsilon^{3-3/m}}+O(1);
\end{align*}
if $m=\frac{5}{3}$,
\begin{align*}
\mathrm{T}=&-\frac{18}{5}\mu\omega_{0}|\ln\varepsilon|-\frac{R\omega_{0}\alpha_{33}^{(\frac{5}{3})}}{\varepsilon^{1/5}}-\frac{R(U_{1}+\omega_{0}R)\alpha_{11}^{(\frac{5}{3})}}{\varepsilon^{2/5}}+\frac{\omega_{0}r^{2}\beta^{(\frac{5}{3})}}{\varepsilon^{6/5}}+O(1);
\end{align*}
if $\frac{5}{3}<m<3$,
\begin{align*}
\mathrm{T}=&-\frac{R(U_{1}+\omega_{0}R)\alpha_{11}^{(m)}}{\varepsilon^{1-1/m}}-\frac{R\omega_{0}\alpha_{33}^{(m)}}{\varepsilon^{2-3/m}}-\frac{\omega_{0}\alpha_{35}^{(m)}}{\varepsilon^{3-5/m}}+\frac{\omega_{0}r^{2}\beta^{(m)}}{\varepsilon^{3-3/m}}+O(1);
\end{align*}
if $m=3$,
\begin{align*}
\mathrm{T}=&\frac{3}{2}\mu\omega_{0}|\ln\varepsilon|-\frac{R(U_{1}+\omega_{0}R)\alpha_{11}^{(3)}}{\varepsilon^{2/3}}-\frac{R\omega_{0}\alpha_{33}^{(3)}}{\varepsilon}-\frac{\omega_{0}\alpha_{35}^{(3)}}{\varepsilon^{4/3}}+\frac{\omega_{0}r^{2}\beta^{(3)}}{\varepsilon^{2}}+O(1);
\end{align*}
if $m>3$,
\begin{align*}
\mathrm{T}=&\frac{\omega_{0}\alpha_{13}^{(m)}}{\varepsilon^{1-3/m}}-\frac{R(U_{1}+\omega_{0}R)\alpha_{11}^{(m)}}{\varepsilon^{1-1/m}}-\frac{R\omega_{0}\alpha_{33}^{(m)}}{\varepsilon^{2-3/m}}-\frac{\omega_{0}\alpha_{35}^{(m)}}{\varepsilon^{3-5/m}}+\frac{\omega_{0}r^{2}\beta^{(m)}}{\varepsilon^{3-3/m}}+O(1),
\end{align*}
where
\begin{align*}
\alpha_{11}^{(m)}&=2\mu\Gamma_{11}^{(m)},\;\alpha_{13}^{(m)}=\frac{\mu}{2m}\left((18+3m)\Gamma_{13}^{(m)}+6m\Gamma_{23}^{(m)}\right),\;\alpha_{33}^{(m)}=3\mu\Gamma_{33}^{(m)},\\
\alpha_{35}^{(m)}&=3\mu\Gamma_{35}^{(m)},\;\beta^{(m)}=3\mu(1+2^{-2}r^{2m-2}-Rr^{m-2})\Gamma_{33}^{(m)},
\end{align*}
and the remainder terms of order $O(\mathbf{1})$ or $O(1)$ depend only on $\mu$, $R$, $r$, $U_{i}$ and $\omega_{0}$, but not on $\varepsilon$.
\end{theorem}

\begin{remark}
Similarly as in $3\mathrm{D}$, Theorem $\ref{thm2}$ shows that the largest blow-up rate $\varepsilon^{3/m-3}$ created by rotational motion appears all in $x_{1}$-axis and $x_{2}$-axis, while the rate $\varepsilon^{3/m-3}$ caused by the translational motion only appears in $x_{2}$-axis. So we claim that the effect of rotational motion on the hydrodynamic force in $x_{1}$-axis is greater than linear motion. Besides, we see that for the $m$-convex particles, the greatest blow-up rate $\varepsilon^{4/m-3}$ in $3\mathrm{D}$ is less than $\varepsilon^{3/m-3}$ in $2\mathrm{D}$.

\end{remark}

\section{Approximation method in $3\mathrm{D}$}
\subsection{Variational formulation}
The solution $\mathbf{u}$ of (\ref{11.3}) minimizes the following functional
\begin{align}\label{5.3}
\mathbf{u}&=\mathrm{arg}\min_{\mathbf{v}\in\mathcal{A}}I_{\Omega}[\mathbf{v}],~~I_{\Omega}[\mathbf{v}]=\frac{1}{2}\int_{\Omega}\boldsymbol{\sigma}(\mathbf{v}):\mathbf{D(v)}\,dx,
\end{align}
where $$\mathcal{A}=\{\mathbf{v}\in\mathbf{H}^{1}(\Omega):~\nabla\cdot\mathbf{v}=0~\mathrm{in}~\Omega,~\mathbf{v}|_{\partial D_{1}}=\mathbf{U}+\boldsymbol{\omega}\times\boldsymbol{\nu},~\mathbf{v}|_{\partial D_{2}}=\mathbf{0}\}.$$
In order to derive (\ref{11.4}) and (\ref{11.5}), it suffices to construct a function in $\mathcal{A}$ to capture the singularity of $\mathbf{F}$ and $\mathbf{T}$.

\subsection{$3\mathrm{D}$ problem in the narrow region between $D_{1}$ and $D_{2}$}

Since the blow-up of the forces only happens in the narrow region between $D_{1}$ and $D_{2},$ we focus on there the asymptotics of the forces. According to (\ref{ZAD01}), we can define the corresponding cylindrical region $\Omega_{r}$ of radius $r$ between the particles $D_{1}$ and $D_{2}$ by
$$\Omega_{r}=\left\{(x',x_{3})\in\mathbb{R}^{3}:~ -\frac{\varepsilon}{2}-\frac{1}{2}|x'|^{m}<x_{3}<\frac{\varepsilon}{2}+\frac{1}{2}|x'|^{m},~ |x'|<r\right\},$$
whose top, bottom and lateral boundaries are given by
\begin{align*}
\Sigma^{+}_{r}&=\left\{(x',x_{3})\in\mathbb{R}^{3}:~ x_{3}=\frac{\varepsilon}{2}+\frac{1}{2}|x'|^{m},~ |x'|<r\right\},\\
\Sigma^{-}_{r}&=\left\{(x',x_{3})\in\mathbb{R}^{3}:~ x_{3}=-\frac{\varepsilon}{2}-\frac{1}{2}|x'|^{m},~ |x'|<r\right\},\\
\partial\Omega_{r}&=\left\{(x',x_{3})\in\mathbb{R}^{3}:~ -\frac{\varepsilon}{2}-\frac{1}{2}|x'|^{m}<x_{3}<\frac{\varepsilon}{2}+\frac{1}{2}|x'|^{m},~ |x'|=r\right\},
\end{align*}
respectively. The vertical distance $h(x')$ between $D_{1}$ and $D_{2}$ passing through $(x',0)$ is given by
\begin{align*}
h(x')&=\varepsilon+|x'|^{m},~\quad |x'|<r.
\end{align*}

We consider the following boundary value problem in the narrow region $\Omega_{r}$:
\begin{align}\label{11.10}
\begin{cases}
(\mathrm{a})~~~\nabla\cdot\boldsymbol{\sigma}\mathbf{(u)}=\mathbf{0}&~~~~~\mathrm{in}~\Omega_{r},\\
(\mathrm{b})~~~~~~~~~~\nabla\cdot\mathbf{u}=0&~~~~~\mathrm{in}~\Omega_{r},\\
(\mathrm{c})~~~~~~~~~~~~~~~~\mathbf{u}=\mathbf{U}+\boldsymbol{\omega}\times\boldsymbol{\nu}&~~~~~\mathrm{on}~\Sigma^{+}_{r},\\
(\mathrm{d})~~~~~~~~~~~~~~~~\mathbf{u}=\mathbf{0}&~~~~~\mathrm{on}~\Sigma^{-}_{r},\\
(\mathrm{e})~\boldsymbol{\sigma}\mathbf{(u)n}_{\partial\Omega_{r}}=\mathbf{0}&~~~~~\mathrm{on}~\partial\Omega_{r},
\end{cases}
\end{align}
where $\mathbf{u}\in\mathbf{H}^{1}(\Omega_{r})$ and $\mathbf{n}_{\partial\Omega_{r}}$ is the unit normal vector to the lateral boundary $\partial\Omega_{r}.$

Next, the major work is to calculate the hydrodynamic force $\mathbf{F}$ and the torque $\mathbf{T}$ on the top boundary $\Sigma_{r}^{+}$ of the narrow region $\Omega_{r}$ as follows:
\begin{align}
\mathbf{F}&=\int_{\Sigma^{+}_{r}}\boldsymbol{\sigma(u)n}\,dS,\quad\mathbf{T}=\int_{\Sigma^{+}_{r}}\boldsymbol{\nu}\times\boldsymbol{\sigma(u)n}\,dS.\label{A11.5}
\end{align}

\subsection{Decomposition of velocity}

Note that the distance $\varepsilon$ between particles $D_{1}$ and $D_{2}$ is much smaller than the radius $r$ of the disk in the plane $x_{3}=0.$ Then we can apply lubrication theory to construct an approximation for the solution $\mathbf{u}$ of $(\ref{11.10})$.

We write
\begin{equation}\label{11.12}
\boldsymbol{\nu}=\mathbf{x}-\mathbf{x}_{D_{1}}=x_{1}\mathbf{e}_{1}+x_{2}\mathbf{e}_{2}+\left(\frac{1}{2}|x'|^{m}-R\right)\mathbf{e}_{3}.
\end{equation}
Denote by $\mathbf{u}|_{\Sigma^{+}_{r}}$ the velocity $\mathbf{u}$ on the top boundary $\Sigma^{+}_{r}$ of $\Omega_{r}$. By $(\ref{11.12})$, we have
\begin{align}\label{11.13}
\mathbf{u}|_{\Sigma^{+}_{r}}=&\mathbf{U}+\boldsymbol{\omega}\times\boldsymbol{\nu}\notag\\
=&\left[U_{1}+\omega_{2}\bigg(\frac{1}{2}|x'|^{m}-R\bigg)-\omega_{3}x_{2}\right]\mathbf{e}_{1}+\left[U_{2}-\omega_{1}\bigg(\frac{1}{2}|x'|^{m}-R\bigg)+\omega_{3}x_{1}\right]\mathbf{e}_{2}\notag\\
&+\left(U_{3}+\omega_{1}x_{2}-\omega_{2}x_{1}\right)\mathbf{e}_{3}.
\end{align}
According to $(\ref{11.13})$ and $\mathbf{u}|_{\Sigma^{-}_{r}}=\mathbf{0}$, we rearrange the terms above as follows:
\begin{align}
\mathbf{u}|_{\Sigma^{+}_{r}}=&\frac{1}{2}(\mathbf{U}+\boldsymbol{\omega}\times\boldsymbol{\nu})+\frac{1}{2}(U_{1}-\omega_{2}R)\,\mathbf{e}_{1}+\frac{1}{2}(U_{2}+\omega_{1}R)\,\mathbf{e}_{2}+\frac{1}{2}U_{3}\mathbf{e}_{3}\notag\\
&+\frac{1}{2}(\omega_{1}x_{2}-\omega_{2}x_{1})\,\mathbf{e}_{3}+\frac{1}{2}\left(\frac{1}{2}\omega_{2}|x'|^{m}-\omega_{3}x_{2}\right)\mathbf{e}_{1}\notag\\
&+\frac{1}{2}\left(-\frac{1}{2}\omega_{1}|x'|^{m}+\omega_{3}x_{1}\right)\mathbf{e}_{2},\label{11.14}\\
\mathbf{u}|_{\Sigma^{-}_{r}}=&\frac{1}{2}(\mathbf{U}+\boldsymbol{\omega}\times\boldsymbol{\nu})-\frac{1}{2}(U_{1}-\omega_{2}R)\,\mathbf{e}_{1}-\frac{1}{2}(U_{2}+\omega_{1}R)\,\mathbf{e}_{2}-\frac{1}{2}U_{3}\mathbf{e}_{3}\notag\\
&-\frac{1}{2}(\omega_{1}x_{2}-\omega_{2}x_{1})\,\mathbf{e}_{3}-\frac{1}{2}\left(\frac{1}{2}\omega_{2}|x'|^{m}-\omega_{3}x_{2}\right)\mathbf{e}_{1}\notag\\
&-\frac{1}{2}\left(-\frac{1}{2}\omega_{1}|x'|^{m}+\omega_{3}x_{1}\right)\mathbf{e}_{2}.\label{C11.14}
\end{align}
Similarly as in \cite{G2016}, making use of the linearity of the Stokes flow and (\ref{11.14})--(\ref{C11.14}), we split the velocity $\mathbf{u}$ as follows:
\begin{equation}\label{11.15}
\mathbf{u}=\sum\limits^{6}_{i=0}\mathbf{u}^{(i)},\quad\mathrm{in}\;\Omega_{r},
\end{equation}
where $\mathbf{u}^{(i)}$, $i\in\{0,\cdots,6\}$ solve problem $((\ref{11.10})\mathrm{(a)}),~((\ref{11.10})\mathrm{(b)}),~((\ref{11.10})(\mathrm{e}))$ and satisfy
\begin{align}\label{11.16}
&(0)~\mathbf{u}^{(0)}|_{\Sigma^{\pm}_{r}}=\frac{1}{2}\left(\mathbf{U}+\mathbf{\omega}\times\boldsymbol{\nu}\right),\notag\\
&(1)~\mathbf{u}^{(1)}|_{\Sigma^{\pm}_{r}}=\pm\frac{1}{2}(U_{1}-\omega_{2}R)\,\mathbf{e}_{1},\notag\\
&(2)~\mathbf{u}^{(2)}|_{\Sigma^{\pm}_{r}}=\pm\frac{1}{2}(U_{2}+\omega_{1}R)\,\mathbf{e}_{2},\notag\\
&(3)~\mathbf{u}^{(3)}|_{\Sigma^{\pm}_{r}}=\pm\frac{1}{2}U_{3}\mathbf{e}_{3},\\
&(4)~\mathbf{u}^{(4)}|_{\Sigma^{\pm}_{r}}=\pm\frac{1}{2}\left(-\omega_{3}x_{2}\mathbf{e}_{1}+\omega_{3}x_{1}\mathbf{e}_{2}\right),\notag\\
&(5)~\mathbf{u}^{(5)}|_{\Sigma^{\pm}_{r}}=\pm\frac{1}{4}|x'|^{m}\left(\omega_{2}\mathbf{e}_{1}-\omega_{1}\mathbf{e}_{2}\right),\notag\\
&(6)~\mathbf{u}^{(6)}|_{\Sigma^{\pm}_{r}}=\pm\frac{1}{2}\left(\omega_{1}x_{2}-\omega_{2}x_{1}\right)\mathbf{e}_{3}.\notag
\end{align}
Then utilizing linearity of the forces in $(\ref{A11.5})$, we have
\begin{align}\label{11.17}
\begin{cases}
\mathbf{F}=\sum\limits^{6}_{i=0}\int_{\Sigma^{+}_{r}}\boldsymbol{\sigma}(\mathbf{u}^{(i)})\mathbf{n}\,dS=:\sum\limits^{6}_{i=0}\mathbf{F}^{(i)},\\
\mathbf{T}=\sum\limits^{6}_{i=0}\int_{\Sigma^{+}_{r}}\boldsymbol{\nu}\times\boldsymbol{\sigma}(\mathbf{u}^{(i)})\mathbf{n}\,dS=:\sum\limits^{6}_{i=0}\mathbf{T}^{(i)}.
\end{cases}
\end{align}

Now we explain the implications of the terms in the decomposition $(\ref{11.16})$. The flow velocity $\mathbf{u}^{(1)}$ is induced by the shear motion of the fluid between $\Sigma^{+}_{r}$ and $\Sigma^{-}_{r}$ along the $x_{1}$-axis with curved surfaces velocities $\pm\frac{1}{2}(U_{1}-\omega_{2}R)\mathbf{e}_{1}$, while $\mathbf{u}^{(2)}$ is created by the shear motion of $\Sigma^{+}_{r}$ and $\Sigma^{-}_{r}$ along the $x_{2}$-axis with velocities $\pm\frac{1}{2}(U_{2}+\omega_{1}R)\mathbf{e}_{2}$. $\mathbf{u}^{(3)}$ and $\mathbf{u}^{(6)}$ are the fluid velocities due to the squeeze motion of $\Sigma^{+}_{r}$ and $\Sigma^{-}_{r}$ with velocities $\pm\frac{1}{2}U_{3}$ and their rotational motion with velocities $\pm\frac{1}{2}\left(\omega_{1}x_{2}-\omega_{2}x_{1}\right)\mathbf{e}_{3}$, respectively. In the remaining terms, $\mathbf{u}^{(0)}$, $\mathbf{u}^{(4)}$ and $\mathbf{u}^{(5)}$ represent flows caused by parallel motion of particles or rotation.

In order to construct the lubrication approximation for the solution $\mathbf{u}$ of $(\ref{11.10})$, we utilize decomposition $(\ref{11.16})$ to construct the approximations for all $\mathbf{u}^{(i)},\,i=0,\cdots,6$. These solutions $\mathbf{u}^{(i)}$ of problems $((\ref{11.10})(a))$, $((\ref{11.10})(b))$, $((\ref{11.10})(e))$ and $(\ref{11.16})$ minimize the corresponding energy functional given by
\begin{align}\label{5.1}
\mathbf{u}^{(i)}&=\mathrm{arg}\min_{\mathbf{v}\in\mathcal{A}_{i}}\mu\int_{\Omega_{r}}\mathbf{D(v)}:\mathbf{D(v)}~dx,
\end{align}
where
\begin{align*}
\mathcal{A}_{i}&=\left\{\mathbf{v}\in\mathbf{H}^{1}(\Omega_{r}):\,\nabla\cdot\mathbf{v}=\mathbf{0}~\mathrm{in}~\Omega_{r},\,\mathbf{v}|_{\Sigma_{r}^{\pm}}=\mathbf{u}^{(i)}~\mathrm{with}~\mathbf{u}^{(i)}~\mathrm{given~by}~(\text{(\ref{11.16})(i)})\right\}.
\end{align*}

We can select an approximation $\overline{\mathbf{u}}^{(i)}$ for the velocity $\mathbf{u}^{(i)}$ from the set $\mathcal{A}_{i}$ and assume its components are polynomials with respect to the variable $x_{3}$.
Combining $((\ref{11.10})(b))$ and $(\ref{11.16})$, we can determine the coefficients of these polynomials. By using
\begin{align}\label{10.1}
\nabla p^{(i)}&=2\mu\nabla\cdot\mathbf{D}(\mathbf{u}^{(i)}),~~~\quad i=0,\cdots,6,
\end{align}
we can recover an approximation $\overline{p}^{(i)}$ for $p^{(i)}$. Moreover, we write the approximation for the hydrodynamic force $\mathbf{F}^{(i)}$ and the torque $\mathbf{T}^{(i)}$ as $\overline{\mathbf{F}}^{(i)}$ and $\overline{\mathbf{T}}^{(i)}$, respectively. That is, in the following we need to calculate
\begin{align}\label{hydro11.179}
\begin{cases}
\overline{\mathbf{F}}=\sum\limits^{6}_{i=0}\int_{\Sigma^{+}_{r}}\boldsymbol{\sigma}(\overline{\mathbf{u}}^{(i)})\mathbf{n}\,dS=:\sum\limits^{6}_{i=0}\overline{\mathbf{F}}^{(i)},\\
\overline{\mathbf{T}}=\sum\limits^{6}_{i=0}\int_{\Sigma^{+}_{r}}\boldsymbol{\nu}\times\boldsymbol{\sigma}(\overline{\mathbf{u}}^{(i)})\mathbf{n}\,dS=:\sum\limits^{6}_{i=0}\overline{\mathbf{T}}^{(i)}.
\end{cases}
\end{align}

\section{Proof of Theorem \ref{thm1}}
In the following, we denote the definite integrals as follows: for $r>0$,
\begin{align*}
\Phi_{ij}^{(m)}(r;\varepsilon)=&\int_{0}^{r}\frac{t^{j}}{(\varepsilon+t^{m})^{i}}\,dt,\quad\mathrm{for}\;i>0,\,j\geq0.
\end{align*}
For $i=0,1,2,...,6$, we denote
\begin{align*}
\overline{\mathbf{F}}^{(i)}=(\overline{F}^{(i)}_{1},\overline{F}^{(i)}_{2},\overline{F}^{(i)}_{3}),\quad\overline{\mathbf{T}}^{(i)}=(\overline{T}^{(i)}_{1},\overline{T}^{(i)}_{2},\overline{T}^{(i)}_{3}).
\end{align*}
We use $\partial_{i}:=\partial_{x_{i}}$, $i=1,2,3$.

The proof of Theorem \ref{thm1} is divided into eight subsections in the following.

\subsection{Asymptotics of $\overline{\mathbf{F}}^{(1)}$ and $\overline{\mathbf{T}}^{(1)}$}
It follows from condition ((\ref{11.16})(1)) that the first component of velocity $\mathbf{u}^{(1)}$ should be an odd function with respect to $x_{3}$. This implies that we can construct an approximation for $\mathbf{u}^{(1)}$ in the following form:
\begin{align}\label{11.18}
\overline{\mathbf{u}}^{(1)}&=(U_{1}-\omega_{2}R)
\left(
\begin{array}{c}
          H(x')x_{3} \\
          0       \\
          -A(x')-B(x')\dfrac{x_{3}^{2}}{2}
\end{array}
 \right).
\end{align}
Due to conditions $((\ref{11.10})(\mathrm{b}))$ and ((\ref{11.16})(1)), we have
\begin{align}\label{11.19}
H(x')=\frac{1}{h(x')},\quad B(x')=\partial_{1}H(x'),\quad A(x')=-\frac{1}{8}B(x')h^{2}(x').
\end{align}
Here we choose any constant $C$ as the approximation for pressure $p^{(1)}$, that is, $\overline{p}^{(1)}=C.$

A direct calculation gives that
\begin{align*}
\nabla\overline{\mathbf{u}}^{(1)}&=(U_{1}-\omega_{2}R)
\left(
\begin{array}{ccc}
           (\partial_{1}H) x_{3} &(\partial_{2}H) x_{3} & H\\
           0 & 0 & 0\\
           -\partial_{1}A-(\partial_{1}B)\dfrac{x_{3}^{2}}{2} & -\partial_{2}A-(\partial_{2}B)\dfrac{x_{3}^{2}}{2} & -Bx_{3}
\end{array}
 \right),\\
D\big(\overline{\mathbf{u}}^{(1)}\big)&=\frac{U_{1}-\omega_{2}R}{2}
\left(
\begin{array}{ccc}
           2(\partial_{1}H)x_{3} &(\partial_{2}H) x_{3} & H-\partial_{1}A-(\partial_{1}B)\dfrac{x_{3}^{2}}{2}\\
           (\partial_{2}H) x_{3} & 0 & -\partial_{2}A-(\partial_{2}B)\dfrac{x_{3}^{2}}{2}\\
           H-\partial_{1}A-(\partial_{1}B)\dfrac{x_{3}^{2}}{2} & -\partial_{2}A-(\partial_{2}B)\dfrac{x_{3}^{2}}{2} & -2Bx_{3}
\end{array}
 \right).
\end{align*}
Thus
\begin{align}
\boldsymbol{\sigma}\big(\overline{\mathbf{u}}^{(1)}\big)\mathbf{n}:=&\mu(U_{1}-\omega_{2}R)\left(I_{1}^{(1)}\mathbf{e}_{1}+I_{2}^{(1)}\mathbf{e}_{2}+I_{3}^{(1)}\mathbf{e}_{3}\right),\label{11.20}\\
\boldsymbol{\nu}\times\boldsymbol{\sigma}\big(\overline{\mathbf{u}}^{(1)}\big)\mathbf{n}=&\mu(U_{1}-\omega_{2}R)
\left|
\begin{array}{ccc}
           \mathbf{e}_{1} & \mathbf{e}_{2} & \mathbf{e}_{3}\\
           \nu_{1} & \nu_{2} & \nu_{3}\\
           I_{1}^{(1)} & I_{2}^{(1)} & I_{3}^{(1)}
\end{array}
 \right|\notag\\
:=&\mu(U_{1}-\omega_{2}R)\left(J_{1}^{(1)}\mathbf{e}_{1}+J_{2}^{(1)}\mathbf{e}_{2}+J_{3}^{(1)}\mathbf{e}_{3}\right),\label{11.21}
\end{align}
where
\begin{align*}
I_{1}^{(1)}&=2(\partial_{1}H)x_{3}n_{1}-\frac{\overline{p}^{(1)}}{\mu(U_{1}-\omega_{2}R)}n_{1}+(\partial_{2}H)x_{3}n_{2}+\bigg(H-\partial_{1}A-(\partial_{1}B)\frac{x_{3}^{2}}{2}\bigg)n_{3},\\
I_{2}^{(1)}&=(\partial_{2}H)x_{3}n_{1}-\frac{\overline{p}^{(1)}}{\mu(U_{1}-\omega_{2}R)}n_{2}+\left(-\partial_{2}A-(\partial_{2}B)\frac{x_{3}^{2}}{2}\right)n_{3},\\
I_{3}^{(1)}&=\left(H-\partial_{1}A-(\partial_{1}B)\frac{x_{3}^{2}}{2}\right)n_{1}+\left(-\partial_{2}A-(\partial_{2}B)\frac{x_{3}^{2}}{2}\right)n_{2}-2Bx_{3}n_{3}-\frac{\overline{p}^{(1)}}{\mu(U_{1}-\omega_{2}R)}n_{3},\\
J_{1}^{(1)}&=\nu_{2}I_{3}^{(1)}-\nu_{3}I_{2}^{(1)},~~J_{2}^{(1)}=\nu_{3}I_{1}^{(1)}-\nu_{1}I_{3}^{(1)},~~J_{3}=\nu_{1}I_{2}^{(1)}-\nu_{2}I_{1}^{(1)}.
\end{align*}
In view of \eqref{hydro11.179}, our objective in the following is to calculate the following surface integrals: for $i=1,2,3$,
\begin{align}\label{GBD01}
\overline{F}^{(1)}_{i}=\mu(U_{1}-\omega_{2}R)\int_{\Sigma^{+}_{r}}I_{i}^{(1)}\,dS,\quad \overline{T}^{(1)}_{i}=\mu(U_{1}-\omega_{2}R)\int_{\Sigma^{+}_{r}}J_{i}^{(1)}\,dS.
\end{align}

First, observe that the unit normal $\mathbf{n}$ to $\Sigma^{+}_{r}$ is given by
$$\mathbf{n}=n_{1}\mathbf{e}_{1}+n_{2}\mathbf{e}_{2}+n_{3}\mathbf{e}_{3},$$
where
\begin{align}\label{11.22}
\begin{cases}
n_{1}=\dfrac{m}{2}\dfrac{x_{1}|x'|^{m-2}}{\sqrt{1+\frac{m^{2}}{4}|x'|^{2m-2}}},\\
n_{2}=\dfrac{m}{2}\dfrac{x_{2}|x'|^{m-2}}{\sqrt{1+\frac{m^{2}}{4}|x'|^{2m-2}}},\\
n_{3}=-\dfrac{1}{\sqrt{1+\frac{m^{2}}{4}|x'|^{2m-2}}}.
\end{cases}
\end{align}
Take the calculation of $\overline{F}^{(1)}_{1}$ defined in \eqref{GBD01} for example. On the top boundary $\Sigma^{+}_{r}$ of $\Omega_{r}$, we first split $I_{1}^{(1)}$ into three parts as follows:
\begin{align*}
I_{1}^{(1)}=I_{11}^{(1)}+I_{12}^{(1)}+I_{13}^{(1)},
\end{align*}
where
\begin{align*}
I_{11}^{(1)}=&2(\partial_{1}H)x_{3}n_{1}+(\partial_{2}H)x_{3}n_{2},\quad I_{12}^{(1)}=-\frac{\overline{p}^{(1)}}{\mu(U_{1}-\omega_{2}R)}n_{1},\notag\\
I_{13}^{(1)}=&\bigg(H-\partial_{1}A-(\partial_{1}B)\frac{x_{3}^{2}}{2}\bigg)n_{3}.
\end{align*}
Then, it follows from (\ref{11.19}) and \eqref{11.22} that
\begin{align*}
\int_{\Sigma^{+}_{r}} I_{11}^{(1)}\,dS=&-\frac{m^{2}}{4}\int_{|x'|<r}\frac{(2x_{1}^{2}+x_{2}^{2})|x'|^{2m-4}}{\varepsilon+|x'|^{m}}\,dx'=O(1),\\
\int_{\Sigma^{+}_{r}} I_{12}^{(1)}\,dS=&-\frac{m\overline{p}^{(1)}}{2\mu(U_{1}-\omega_{2}R)}\int_{|x'|<r}x_{1}|x'|^{m-2}dx'=0,
\end{align*}
while
\begin{align}\label{11.25}
\int_{\Sigma^{+}_{r}}I_{13}^{(1)}\,dS
=&-\int_{\{|x'|<r\}}\left(\frac{1}{\varepsilon+|x'|^{m}}-\frac{m^{2}}{4}\frac{x_{1}^{2}|x'|^{2m-4}}{\varepsilon+|x'|^{m}}\right)dx'\notag\\
=&-2\pi \Phi_{11}^{(m)}(r;\varepsilon)+O(1).
\end{align}
Then
\begin{align}\label{11.25}
\overline{F}^{(1)}_{1}=\mu(U_{1}-\omega_{2}R)\int_{\Sigma^{+}_{r}}I_{1}^{(1)}\,dS=-2\pi\mu(U_{1}-\omega_{2}R) \Phi_{11}^{(m)}(r;\varepsilon)+O(1).
\end{align}
We would like to remark that the integral of $I_{12}^{(1)}$ becomes zero due to the parity of integrand and the symmetry of domain. This fact also makes $\overline{\mathbf{F}}^{(1)}$ in other directions $e_{2}$ and $e_{3}$ no singularities. That is, the singularities of $\overline{\mathbf{F}}^{(1)}$ only appear in the direction of $\mathbf{e}_{1}$. By the same argument, the singularities of $\overline{\mathbf{T}}^{(1)}$ only appear in the direction of $\mathbf{e}_{2}$,
\begin{align}\label{11.26}
\overline{T}^{(1)}_{2}=&\mu(U_{1}-\omega_{2}R)\int_{\Sigma^{+}_{r}}\bigg(H-\partial_{1}A-(\partial_{1}B)\frac{x_{3}^{2}}{2}\bigg)n_{3}\nu_{3}\,dS+O(1)\notag\\
&=2\pi\mu R(U_{1}-\omega_{2}R)\Phi_{11}^{(m)}(r;\varepsilon)+O(1).
\end{align}
Therefore, combining (\ref{11.20})--(\ref{11.21}) and (\ref{11.25})--(\ref{11.26}), we obtain that for $m=2$,
\begin{align}
\overline{\mathbf{F}}^{(1)}&=-\pi\mu(U_{1}-\omega_{2}R)|\ln\varepsilon|\,\mathbf{e}_{1}+O(\mathbf{1}),\label{21.010}\\
\overline{\mathbf{T}}^{(1)}&=\pi\mu R(U_{1}-\omega_{2}R)|\ln\varepsilon|\,\mathbf{e}_{2}+O(\mathbf{1});\label{21.011}
\end{align}
for $m>2$,
\begin{align}
\overline{\mathbf{F}}^{(1)}&=-\frac{2\pi\mu(U_{1}-\omega_{2}R)\Gamma_{12}^{(m)}}{\varepsilon^{1-{2}/{m}}}\,\mathbf{e}_{1}+O(\mathbf{1}),\label{11.27}\\
\overline{\mathbf{T}}^{(1)}&=\frac{2\pi\mu R(U_{1}-\omega_{2}R)\Gamma_{12}^{(m)}}{\varepsilon^{1-{2}/{m}}}\,\mathbf{e}_{2}+O(\mathbf{1}).\label{11.28}
\end{align}

Since the method of calculation is the same for all $\overline{\mathbf{F}}^{(i)}$ and $\overline{\mathbf{T}}^{(i)}$, $i=0,1,2,...,6$, we just give the calculations of the singular terms in the following sections for the sake of simplicity.
\subsection{Asymptotics of $\overline{\mathbf{F}}^{(2)}$ and $\overline{\mathbf{T}}^{(2)}$}\label{JNK}

Similarly as before, we assume
\begin{align}\label{11.29}
\overline{\mathbf{u}}^{(2)}&=(U_{2}+\omega_{1}R)
\left(
\begin{array}{c}
             0    \\
          H(x')x_{3} \\
          -A(x')-B(x')\dfrac{x_{3}^{2}}{2}
\end{array}
 \right).
\end{align}
Utilizing $((\ref{11.10})(\mathrm{b}))$ and $((\ref{11.16})(2))$, we have
\begin{align}\label{11.30}
H(x')=\frac{1}{h(x')},\quad B(x')=\partial_{1}H(x'),\quad A(x')=-\frac{1}{8}B(x')h^{2}(x').
\end{align}
By the same argument as in Section \ref{JNK}, letting $\overline{p}^{(2)}=C,$ we deduce that for $m=2$,
\begin{align}
\overline{\mathbf{F}}^{(2)}&=-\pi\mu(U_{2}+\omega_{1}R)|\ln\varepsilon|\,\mathbf{e}_{2}+O(\mathbf{1}),\label{21.012}\\
\overline{\mathbf{T}}^{(2)}&=-\pi\mu R(U_{2}+\omega_{1}R)|\ln\varepsilon|\,\mathbf{e}_{1}+O(\mathbf{1});\label{21.013}
\end{align}
for $m>2$,
\begin{align}
\overline{\mathbf{F}}^{(2)}&=-\frac{2\pi\mu(U_{2}+\omega_{1}R)\Gamma_{12}^{(m)}}{\varepsilon^{1-{2}/{m}}}\,\mathbf{e}_{2}+O(\mathbf{1}),\label{11.31}\\
\overline{\mathbf{T}}^{(2)}&=-\frac{2\pi\mu R(U_{2}+\omega_{1}R)\Gamma_{12}^{(m)}}{\varepsilon^{1-{2}/{m}}}\,\mathbf{e}_{1}+O(\mathbf{1}).\label{11.32}
\end{align}

\subsection{Asymptotics of $\overline{\mathbf{F}}^{(3)}$ and $\overline{\mathbf{T}}^{(3)}$}

Choose
\begin{align}\label{11.33}
\overline{\mathbf{u}}^{(3)}&=U_{3}
\left(
\begin{array}{c}
          -A_{1}(x')-3B_{1}(x')x_{3}^{2}    \\\\
          -A_{2}(x')-3B_{2}(x')x_{3}^{2} \\\\
          A_{3}(x')x_{3}+B_{3}(x')x_{3}^{3}
\end{array}
 \right).
\end{align}
It follows from $((\ref{11.10})(\mathrm{b}))$ and $((\ref{11.16})(3))$ that
\begin{align}
A_{1}&=\frac{3}{4}\frac{x_{1}}{h},\quad A_{2}=\frac{3}{4}\frac{x_{2}}{h},\quad B_{1}=-\frac{x_{1}}{h^{3}},\quad B_{2}=-\frac{x_{2}}{h^{3}},\label{11.34}\\
A_{3}&=\partial_{1}A_{1}+\partial_{2}A_{2},\quad B_{3}=\partial_{1}B_{1}+\partial_{2}B_{2}.\label{11.344}
\end{align}

In view of $(\ref{10.1})$, instead of choosing the constant pressure as before, we let
\begin{equation}\label{11.35}
\overline{p}^{(3)}=\mu U_{3}\left(-A_{3}(x')+3B_{3}(x')x_{3}^{2}-6G(x')\right),
\end{equation}
where
\begin{equation}\label{11.38}
G(x')=-\frac{1}{2}\int^{x_{1}^{2}+x_{2}^{2}}_{r^{2}}\frac{1}{(\varepsilon+s^{\frac{m}{2}})^{3}}\,ds,
\end{equation}
satisfying $\partial_{1}G=B_{1}$ and $\partial_{2}G=B_{2}.$

By definition,
\begin{align*}
\nabla\overline{\mathbf{u}}^{(3)}&=U_{3}
\left(
\begin{array}{ccc}
           -\partial_{1}A_{1}-3(\partial_{1}B_{1})x_{3}^{2} & -\partial_{2}A_{1}-3(\partial_{2}B_{1})x_{3}^{2} & -6B_{1}x_{3}\\\\
           -\partial_{1}A_{2}-3(\partial_{1}B_{2})x_{3}^{2} & -\partial_{2}A_{2}-3(\partial_{2}B_{2})x_{3}^{2} & -6B_{2}x_{3}\\\\
           (\partial_{1}A_{3})x_{3}+(\partial_{1}B_{3})x_{3}^{3} & (\partial_{2}A_{3})x_{3}+(\partial_{2}B_{3})x_{3}^{3} & A_{3}+3B_{3}x_{3}^{2}
\end{array}
 \right).
\end{align*}
Consequently,
\begin{align}
\boldsymbol{\sigma}\big(\overline{\mathbf{u}}^{(3)}\big)\mathbf{n}:=&\mu U_{3}\left(I_{1}^{(3)}\mathbf{e}_{1}+I_{2}^{(3)}\mathbf{e}_{2}+I_{3}^{(3)}\mathbf{e}_{3}\right),\label{11.36}\\
\boldsymbol{\nu}\times\boldsymbol{\sigma}\big(\overline{\mathbf{u}}^{(3)}\big)\mathbf{n}=&\mu U_{3}
\left|
\begin{array}{ccc}
           \mathbf{e}_{1} & \mathbf{e}_{2} & \mathbf{e}_{3}\\
           \nu_{1} & \nu_{2} & \nu_{3}\\
           I_{1}^{(3)} & I_{2}^{(3)} & I_{3}^{(3)}
\end{array}
 \right|\notag\\
:=&\mu U_{3}\left(J_{1}^{(3)}\mathbf{e}_{1}+J_{2}^{(3)}\mathbf{e}_{2}+J_{3}^{(3)}\mathbf{e}_{3}\right),\label{11.37}
\end{align}
where
\begin{align*}
I_{1}^{(3)}=&2\left(-\partial_{1}A_{1}-3(\partial_{1}B_{1})x_{3}^{2}\right)n_{1}-\frac{\overline{p}^{(3)}}{\mu U_{3}}n_{1}-\left[\partial_{2}A_{1}+\partial_{1}A_{2}+3x_{3}^{2}(\partial_{2}B_{1}+\partial_{1}B_{2})\right]n_{2}\\
&+\left(\partial_{1}A_{3}x_{3}+(\partial_{1}B_{3})x_{3}^{3}-6B_{1}x_{3}\right)n_{3},\\
I_{2}^{(3)}=&-\left[\partial_{2}A_{1}+\partial_{1}A_{2}+3x_{3}^{2}(\partial_{2}B_{1}+\partial_{1}B_{2})\right]n_{1}+2\left(-\partial_{2}A_{2}-3(\partial_{2}B_{2})x_{3}^{2}\right)n_{2}\\
&-\frac{\overline{p}^{(3)}}{\mu U_{3}}n_{2}+\left(\partial_{2}A_{3}x_{3}+\partial_{2}B_{3}x_{3}^{3}-6B_{2}x_{3}\right)n_{3},\\
I_{3}^{(3)}=&\left((\partial_{1}A_{3})x_{3}+(\partial_{1}B_{3})x_{3}^{3}-6B_{1}x_{3}\right)n_{1}+\left((\partial_{2}A_{3})x_{3}+(\partial_{2}B_{3})x_{3}^{3}-6B_{2}x_{3}\right)n_{2}\\
&+2(A_{3}+3B_{3}x_{3}^{2})n_{3}-\frac{\overline{p}^{(3)}}{\mu U_{3}}n_{3},\\
J_{1}^{(3)}=&\nu_{2}I_{3}^{(3)}-\nu_{3}I_{2}^{(3)},~~J_{2}^{(3)}=\nu_{3}I_{1}^{(3)}-\nu_{1}I_{3}^{(3)},~~J_{3}=\nu_{1}I_{2}^{(3)}-\nu_{2}I_{1}^{(3)}.
\end{align*}
By definition, our goal is to calculate the surface integrals as follows: for $i=1,2,3$,
\begin{align}\label{GBD03}
\overline{F}^{(3)}_{i}=\mu U_{3}\int_{\Sigma^{+}_{r}}I_{i}^{(3)}\,dS,\quad \overline{T}^{(3)}_{i}=\mu U_{3}\int_{\Sigma^{+}_{r}}J_{i}^{(3)}\,dS.
\end{align}

In view of (\ref{11.34})--(\ref{11.38}), on the top boundary $\Sigma^{+}_{r}$ we have
\begin{align}
(\partial_{1}A_{3})x_{3}+(\partial_{1}B_{3})x_{3}^{3}-6B_{1}x_{3}&=-\frac{3}{4}\frac{(x_{1}\partial_{1}h+x_{2}\partial_{2}h)\partial_{1}h}{h^{2}}+\frac{3x_{1}}{h^{2}},\label{11.39}\\
(\partial_{2}A_{3})x_{3}+(\partial_{2}B_{3})x_{3}^{3}-6B_{2}x_{3}&=-\frac{3}{4}\frac{(x_{1}\partial_{1}h+x_{2}\partial_{2}h)\partial_{2}h}{h^{2}}+\frac{3x_{2}}{h^{2}},\label{11.40}\\
A_{3}+3B_{3}x_{3}^{2}&=\frac{3}{2}\frac{x_{1}\partial_{1}h+x_{2}\partial_{2}h}{h^{2}},\label{11.41}\\
\frac{\overline{p}^{(3)}}{\mu U_{3}}&=\frac{3(x_{1}\partial_{1}h+x_{2}\partial_{2}h)}{h^{2}}-\frac{1}{h}-6G.\label{11.42}
\end{align}

Similarly as before, by utilizing the parity of integrand and the symmetry of domain, we deduce that the singularity of $\overline{\mathbf{F}}^{(3)}$ only lies in $\mathbf{e}_{3}$, while $\overline{\mathbf{T}}^{(3)}$ has no singularity. By a direct calculation, it follows from (\ref{11.36})--(\ref{11.42}) and integration by parts that
\begin{align}
\int_{\Sigma^{+}_{r}}\left((\partial_{1}A_{3})x_{3}+(\partial_{1}B_{3})x_{3}^{3}-6B_{1}x_{3}\right)n_{1}\,dS=&\frac{3}{2}m\int_{|x'|<r}\frac{x_{1}^{2}|x'|^{m-2}}{(\varepsilon+|x'|^{m})^{2}}dx'+O(1)\notag\\
=&\frac{3}{2}m\pi\int_{0}^{r}\frac{t^{m+1}}{(\varepsilon+t^{m})^{2}}dt+O(1)\notag\\
=&3\pi \Phi_{11}^{(m)}(r;\varepsilon)+O(1),\label{11.43}
\end{align}
and
\begin{align}
&\int_{\Sigma^{+}_{r}}\left((\partial_{2}A_{3})x_{3}+(\partial_{2}B_{3})x_{3}^{3}-6B_{2}x_{3}\right)n_{2}\,dS=3\pi \Phi_{11}^{(m)}(r;\varepsilon)+O(1),\label{11.44}\\
&\int_{\Sigma^{+}_{r}}2(A_{3}+3B_{3}x_{3}^{2})n_{3}\,dS=-12\pi \Phi_{11}^{(m)}(r;\varepsilon)+O(1),\label{11.45}\\
&\int_{\Sigma^{+}_{r}}-\frac{\overline{p}^{(3)}}{\mu U_{3}}n_{3}\,dS=6\pi \Phi_{11}^{(m)}(r;\varepsilon)-\frac{3\pi\Gamma_{34}^{(m)}}{\varepsilon^{3-{4}/{m}}}+O(1).\label{11.46}
\end{align}

According to \eqref{GBD03}, it follows from (\ref{11.43})--(\ref{11.46}) that for $m\geq2$,
\begin{align*}
\overline{F}^{(3)}_{3}=\mu U_{3}\int_{\Sigma^{+}_{r}}I_{3}^{(3)}\,dS=\frac{3\pi\mu U_{3}\Gamma_{34}^{(m)}}{\varepsilon^{3-{4}/{m}}}+O(1).
\end{align*}
Consequently,
\begin{align}
\overline{\mathbf{F}}^{(3)}&=\frac{3\pi\mu U_{3}\Gamma_{34}^{(m)}}{\varepsilon^{3-{4}/{m}}}\,\mathbf{e}_{3}+O(\mathbf{1}),\quad\overline{\mathbf{T}}^{(3)}=\mathbf{0}+O(\mathbf{1}).\label{11.51}
\end{align}

\subsection{Asymptotics of $\overline{\mathbf{F}}^{(4)}$ and $\overline{\mathbf{T}}^{(4)}$}

Assume that
\begin{align}\label{11.52}
\overline{\mathbf{u}}^{(4)}&=\omega_{3}
\left(
\begin{array}{c}
          H_{1}(x')x_{3} \\
          0       \\
          -A_{1}(x')-B_{1}(x')\frac{x_{3}^{2}}{2}
\end{array}
 \right)+
\omega_{3}
\left(
\begin{array}{c}
              0   \\
          H_{2}(x')x_{3} \\
          -A_{2}(x')-B_{2}(x')\frac{x_{3}^{2}}{2}
\end{array}
 \right),
\end{align}
which, in combination with incompressibility and ((\ref{11.16})(4)), yields that
\begin{align}
H_{1}=&-\frac{x_{2}}{h},\quad H_{2}=\frac{x_{1}}{h},\quad B_{1}=\partial_{1}H_{1},\label{11.53}\\
B_{2}=&\partial_{2}H_{2},\quad A_{1}=-\frac{B_{1}h^{2}}{8},\quad A_{2}=-\frac{B_{2}h^{2}}{8}.\label{11.533}
\end{align}
Let $\overline{p}^{(4)}=C$. Similarly as above, a direct calculation gives
\begin{equation}\label{11.54}
\overline{\mathbf{F}}^{(4)}=\mathbf{0}+O(\mathbf{1}),\quad \overline{\mathbf{T}}^{(4)}=\mathbf{0}+O(\mathbf{1}).
\end{equation}

\subsection{Asymptotics of $\overline{\mathbf{F}}^{(5)}$ and $\overline{\mathbf{T}}^{(5)}$}

Let
\begin{align}\label{11.55}
\overline{\mathbf{u}}^{(5)}&=\omega_{2}
\left(
\begin{array}{c}
          H_{1}(x')x_{3} \\
          0       \\
          -A_{1}(x')-B_{1}(x')\dfrac{x_{3}^{2}}{2}
\end{array}
 \right)+
\omega_{1}
\left(
\begin{array}{c}
              0   \\
          H_{2}(x')x_{3} \\
          -A_{2}(x')-B_{2}(x')\dfrac{x_{3}^{2}}{2}
\end{array}
 \right),
\end{align}
where
\begin{align}
A_{1}&=-\frac{\varepsilon}{16}\partial_{1}h,\quad B_{1}=\frac{\varepsilon}{2}\frac{\partial_{1}h}{h^{2}},\quad H_{1}=\frac{1}{2}-\frac{\varepsilon}{2}\frac{1}{h},\label{15.01}\\
A_{2}&=\frac{\varepsilon}{16}\partial_{2}h,\quad B_{2}=-\frac{\varepsilon}{2}\frac{\partial_{2}h}{h^{2}},\quad H_{2}=-\frac{1}{2}+\frac{\varepsilon}{2}\frac{1}{h},\label{15.02}
\end{align}
under conditions $((\ref{11.10})(\mathrm{b}))$ and $((\ref{11.16})(5))$.

Choose $\overline{p}^{(5)}=C.$ Similarly, a direct calculation yields
\begin{align}
\overline{\mathbf{F}}^{(5)}=\mathbf{0}+O(\mathbf{1}),\quad\overline{\mathbf{T}}^{(5)}=\mathbf{0}+O(\mathbf{1}).\label{21.015}
\end{align}
\subsection{Estimates of $\overline{\mathbf{F}}^{(6)}$ and $\overline{\mathbf{T}}^{(6)}$}

As before, we assume
\begin{align}\label{11.60}
\overline{\mathbf{u}}^{(6)}&=
\left(
\begin{array}{c}
          -A_{1}(x')-3B_{1}(x')x_{3}^{2}    \\\\
          -A_{2}(x')-3B_{2}(x')x_{3}^{2} \\\\
          A_{3}(x')x_{3}+B_{3}(x')x_{3}^{3}
\end{array}
 \right).
\end{align}
Making use of $((\ref{11.10})(\mathrm{b}))$ and $((\ref{11.16})(6))$, we deduce that
\begin{equation}\label{11.61}
A_{1}=-\frac{3}{4}\omega_{2}\frac{x_{1}^{2}}{h},\quad A_{2}=\frac{3}{4}\omega_{1}\frac{x_{2}^{2}}{h},\quad B_{1}=\omega_{2}\frac{x_{1}^{2}}{h^{3}},\quad B_{2}=-\omega_{1}\frac{x_{2}^{2}}{h^{3}},
\end{equation}
and $A_{3}$ and $B_{3}$ have the same form as $(\ref{11.344})$.

In light of $(\ref{10.1})$, we choose
\begin{equation}\label{11.62}
\overline{p}^{(6)}=-\mu A_{3}(x')+3\mu B_{3}(x')x_{3}^{2}-6\mu G_{1}(x')-6\mu G_{2}(x')
\end{equation}
with $\partial_{1}G_{1}=B_{1}(x'),~\partial_{2}G_{2}=B_{2}(x')$.

Therefore, by definition, we obtain
\begin{align}
\boldsymbol{\sigma}\big(\overline{\mathbf{u}}^{(6)}\big)\mathbf{n}:=&\mu\left(I_{1}\mathbf{e}_{1}+I_{2}\mathbf{e}_{2}+I_{3}\mathbf{e}_{3}\right),\label{11.63}\\
\boldsymbol{\nu}\times\boldsymbol{\sigma}\big(\overline{\mathbf{u}}^{(6)}\big)\mathbf{n}=&\mu
\left|
\begin{array}{ccc}
           \mathbf{e}_{1} & \mathbf{e}_{2} & \mathbf{e}_{3}\\
           \nu_{1} & \nu_{2} & \nu_{3}\\
           I_{1}^{(6)} & I_{2}^{(6)} & I_{3}^{(6)}
\end{array}
 \right|\notag\\
:=&\mu\left(J_{1}^{(6)}\mathbf{e}_{1}+J_{2}^{(6)}\mathbf{e}_{2}+J_{3}^{(6)}\mathbf{e}_{3}\right),\label{11.64}
\end{align}
where
\begin{align*}
I_{1}^{(6)}=&2\left(-\partial_{1}A_{1}-3(\partial_{1}B_{1})x_{3}^{2}\right)n_{1}-\frac{\overline{p}^{(6)}}{\mu}n_{1}-\left[\partial_{2}A_{1}+\partial_{1}A_{2}+3x_{3}^{2}(\partial_{2} B_{1}+\partial_{1}B_{2})\right]n_{2}\\
&+\left((\partial_{1}A_{3})x_{3}+(\partial_{1}B_{3})x_{3}^{3}-6B_{1}x_{3}\right)n_{3},\\
I_{2}^{(6)}=&-\left[\partial_{2}A_{1}+\partial_{1}A_{2}+3x_{3}^{2}(\partial_{2}B_{1}+\partial_{1}B_{2})\right]n_{1}+2\left(-\partial_{2}A_{2}-3(\partial_{2}B_{2})x_{3}^{2}\right)n_{2}\\
&-\frac{\overline{p}^{(6)}}{\mu}n_{2}+\left((\partial_{2}A_{3})x_{3}+(\partial_{2}B_{3})x_{3}^{3}-6B_{2}x_{3}\right)n_{3},\\
I_{3}^{(6)}=&\left((\partial_{1}A_{3})x_{3}+(\partial_{1}B_{3})x_{3}^{3}-6B_{1}x_{3}\right)n_{1}+\left((\partial_{2}A_{3})x_{3}+(\partial_{2}B_{3})x_{3}^{3}-6B_{2}x_{3}\right)n_{2}\\
&+2(A_{3}+3B_{3}x_{3}^{2})n_{3}-\frac{\overline{p}^{(6)}}{\mu}n_{3},\\
J_{1}^{(6)}=&\nu_{2}I_{3}^{(6)}-\nu_{3}I_{2}^{(6)},~~J_{2}^{(6)}=\nu_{3}I_{1}^{(6)}-\nu_{1}I_{3}^{(6)},~~J_{3}^{(6)}=\nu_{1}I_{2}^{(6)}-\nu_{2}I_{1}^{(6)}.
\end{align*}
By definition, the purpose of this section is to calculate the surface integrals as follows: for $i=1,2,3$,
\begin{align}\label{GBD06}
\overline{F}^{(6)}_{i}=\mu\int_{\Sigma^{+}_{r}}I_{i}^{(6)}\,dS,\quad \overline{T}^{(6)}_{i}=\mu\int_{\Sigma^{+}_{r}}J_{i}^{(6)}\,dS.
\end{align}

Consider
\begin{equation}\label{3.0021}
G_{1}=\omega_{2}\int^{x_{1}}_{r}\frac{t^{2}}{h^{3}(t,x_{2})}\,dt,\quad G_{2}=-\omega_{1}\int^{x_{2}}_{-r}\frac{t^{2}}{h^{3}(x_{1},t)}\,dt.
\end{equation}

On the top boundary $\Sigma^{+}_{r}$, it follows from (\ref{11.61})--(\ref{3.0021}) that
\begin{align}
&(\partial_{1}A_{3})x_{3}+(\partial_{1}B_{3})x_{3}^{3}-6B_{1}x_{3}\notag\\
=&\frac{\omega_{2}}{4}\left(-2+\frac{3x_{1}^{2}(\partial_{1}h)^{2}}{h^{2}}\right)-\frac{3}{4}\omega_{1}\frac{x_{2}^{2}(\partial_{1}h)(\partial_{2}h)}{h^{2}}-3\omega_{2}\frac{x_{1}^{2}}{h^{2}},\label{11.644}\\
&(\partial_{2}A_{3})x_{3}+(\partial_{2}B_{3})x_{3}^{3}-6B_{2}x_{3}\notag\\
=&\frac{3}{4}\omega_{2}\frac{x_{1}^{2}(\partial_{1}h)(\partial_{2}h)}{h^{2}}+\frac{\omega_{1}}{4}\left(2-\frac{3x_{2}^{2}(\partial_{2}h)^{2}}{h^{2}}\right)+3\omega_{1}\frac{x_{2}^{2}}{h^{2}},\label{11.65}
\end{align}
and
\begin{align}
&\frac{\overline{p}^{(6)}}{\mu}=3\omega_{2}\frac{x_{1}h-x_{1}^{2}\partial_{1}h}{h^{2}}-3\omega_{1}\frac{x_{2}h-x_{2}^{2}\partial_{2}h}{h^{2}}-6\omega_{2}\int^{x_{1}}_{r}\frac{t^{2}}{h^{3}(t,x_{2})}dt\notag\\
&\quad\quad\quad+6\omega_{1}\int^{x_{2}}_{-r}\frac{t^{2}}{h^{3}(x_{1},t)}dt.\label{3.008}
\end{align}

We now estimate the singularities of $\overline{F}^{(6)}_{i}$ and $\overline{T}^{(6)}_{i}$, $i=1,2,3$, respectively. For simplicity, $O(1)$ is omitted in the following estimates.

$\mathrm{\mathbf{Step}}~\mathrm{\mathbf{1}}$. Estimates of $\overline{F}^{(6)}_{i},~i=1,2,3$.

First, applying integration by parts, it follows from $(\ref{3.008})$ that
\begin{align}
\int_{\Sigma^{+}_{r}}-\frac{\overline{p}^{(6)}}{3\mu}n_{1}=&\int_{\{|x'|<r\}}\partial_{1}hG_{1}\notag\\
\leq&\omega_{2}\int^{r}_{-r}dx_{2}\int^{r}_{-r}\partial_{1}h\int^{x_{1}}_{r}\frac{t^{2}}{h^{3}(t,x_{2})}dt~dx_{1}\notag\\
\leq&\omega_{2}(\varepsilon+2^{\frac{m}{2}}r^{m})\int^{r}_{-r}\int^{r}_{-r}\frac{x_{1}^{2}}{h^{3}}dx_{1}dx_{2}-\omega_{2}\int^{r}_{-r}\int^{r}_{-r}\frac{x_{1}^{2}}{h^{2}}dx_{1}dx_{2}\notag\\
\leq&\omega_{2}(\varepsilon+2^{\frac{m}{2}}r^{m})\int_{\{|x'|<2r\}}\frac{x_{1}^{2}}{h^{3}}-\omega_{2}\int_{\{|x'|<r\}}\frac{x_{1}^{2}}{h^{2}}\notag\\
\leq&\pi\omega_{2}(\varepsilon+2^{\frac{m}{2}}r^{m})\Phi_{33}^{(m)}(2r;\varepsilon)-\pi\omega_{2}\Phi_{23}^{(m)}(r;\varepsilon).\label{3.009}
\end{align}
Second,
\begin{align}
\int_{\Sigma^{+}_{r}}-\frac{\overline{p}^{(6)}}{3\mu}n_{1}=&\int_{\{|x'|<r\}}\partial_{1}hG_{1}\notag\\
\geq&\omega_{2}\int^{\frac{r}{2}}_{-\frac{r}{2}}dx_{2}\int^{\frac{r}{2}}_{-\frac{r}{2}}\partial_{1}h\int^{x_{1}}_{r}\frac{t^{2}}{h^{3}(t,x_{2})}dt~dx_{1}\notag\\
\geq&\omega_{2}(\varepsilon+2^{-m}r^{m})\int_{\{|x_{1}|,|x_{2}|<\frac{r}{2}\}}\frac{x_{1}^{2}}{h^{3}}-\omega_{2}\int_{\{|x_{1}|,|x_{2}|<\frac{r}{2}\}}\frac{x_{1}^{2}}{h^{2}}\notag\\
\geq&\omega_{2}(\varepsilon+2^{-m}r^{m})\int_{\{|x'|<\frac{r}{2}\}}\frac{x_{1}^{2}}{h^{3}}-\omega_{2}\int_{\{|x'|<r\}}\frac{x_{1}^{2}}{h^{2}}\notag\\
\geq&\pi\omega_{2}(\varepsilon+2^{-m}r^{m})\Phi_{33}^{(m)}(2r;\varepsilon)-\pi\omega_{2}\Phi_{23}^{(m)}(r;\varepsilon).\label{3.0010}
\end{align}
In addition, by using (\ref{11.644}), we obtain
\begin{align}
\int_{\Sigma^{+}_{r}}\left((\partial_{1}A_{3})x_{3}+(\partial_{1}B_{3})x_{3}^{3}-6B_{1}x_{3}\right)n_{3}\,dS
&=3\pi\omega_{2}\Phi_{23}^{(m)}(r;\varepsilon).\label{11.67}
\end{align}
Hence, in light of \eqref{GBD06}, it follows from (\ref{3.009})--(\ref{11.67}) that for $m=2$,
\begin{equation}\label{21.017}
\frac{3\pi\omega_{2}r^{2}}{16\varepsilon}\leq\overline{F}^{(6)}_{1}\leq\frac{3\pi\omega_{2}r^{2}}{2\varepsilon};
\end{equation}
for $m>2$,
\begin{equation}\label{3.0011} \frac{3\pi\mu\omega_{2}\Gamma_{34}^{(m)}}{2}\frac{2^{-m}r^{m}+\varepsilon}{\varepsilon^{3-{4}/{m}}}\leq \overline{F}^{(6)}_{1}\leq\frac{3\pi\mu\omega_{2}\Gamma_{34}^{(m)}}{2}\frac{2^{\frac{m}{2}}r^{m}+\varepsilon}{\varepsilon^{3-{4}/{m}}}.
\end{equation}

Similarly as before, utilizing (\ref{11.65}) and (\ref{3.008}), we deduce that for $m=2$,
\begin{align}
-\frac{3\pi\omega_{1}r^{2}}{2\varepsilon}&\leq\overline{F}^{(6)}_{2}\leq-\frac{3\pi\omega_{1}r^{2}}{16\varepsilon},\label{21.018}\\
\frac{3\pi\mu(\omega_{1}+\omega_{2})r}{4\varepsilon}&\leq\overline{F}^{(6)}_{3}\leq\frac{3\pi\mu(\omega_{1}+\omega_{2})r}{2\varepsilon};\label{21.019}
\end{align}
for $m>2$,
\begin{align} -\frac{3\pi\mu\omega_{1}\Gamma_{34}^{(m)}}{2}\frac{2^{\frac{m}{2}}r^{m}+\varepsilon}{\varepsilon^{3-{4}/{m}}}\leq&\overline{F}^{(6)}_{2}\leq-\frac{3\pi\mu\omega_{1}\Gamma_{34}^{(m)}}{2}\frac{2^{-m}r^{m}+\varepsilon}{\varepsilon^{3-{4}/{m}}},\label{3.0012}\\
\frac{3\pi\mu r(\omega_{1}+\omega_{2})\Gamma_{34}^{(m)}}{2\varepsilon^{3-{4}/{m}}}\leq&\overline{F}^{(6)}_{3}\leq\frac{3\pi\mu r(\omega_{1}+\omega_{2})\Gamma_{34}^{(m)}}{\varepsilon^{3-{4}/{m}}}.\label{3.0013}
\end{align}

$\mathrm{\mathbf{Step}}~\mathrm{\mathbf{2}}$. Estimates of $\overline{T}^{(6)}_{i},\;i=1,2,3$.

Recalling the definition of $\overline{T}^{(6)}_{1}$ and making use of (\ref{11.644})--(\ref{3.008}), we have
\begin{align*}
&\int_{\Sigma^{+}_{r}}\left((\partial_{2}A_{3})x_{3}+(\partial_{2}B_{3})x_{3}^{3}-6B_{2}x_{3}\right)n_{3}\nu_{3}\,dS\\
=&3\pi\omega_{1}R\Phi_{23}^{(m)}(r;\varepsilon)-\frac{6\omega_{1}\pi}{m}\Phi_{13}^{(m)}(r;\varepsilon),\\
&\int_{\Sigma^{+}_{r}}-\left((\partial_{2}A_{3})x_{3}+(\partial_{2}B_{3})x_{3}^{3}-6B_{2}x_{3}\right)n_{2}\nu_{2}=-\frac{9\pi\omega_{1}}{2}\Phi_{13}^{(m)}(r;\varepsilon),\\
&\int_{\Sigma^{+}_{r}}-2(A_{3}+3B_{3}x_{3}^{2})n_{3}\nu_{2}=9\pi\omega_{1}\Phi_{13}^{(m)}(r;\varepsilon),
\end{align*}
and
\begin{align*}
\int_{\Sigma^{+}_{r}}\frac{\overline{p}^{(6)}}{\mu}n_{3}\nu_{2}=&-6\pi\omega_{1}\Phi_{13}^{(m)}(r;\varepsilon)-6\omega_{1}\int_{\{|x'|<r\}}x_{2}\int^{x_{2}}_{-r}\frac{t^{2}}{h^{3}(x_{1},t)}dt,\\
\int_{\Sigma^{+}_{r}}-\frac{\overline{p}^{(6)}}{\mu}n_{2}\nu_{3}=&\frac{6\omega_{1}R+3\omega_{1}\varepsilon}{2}\int_{\{|x'|<r\}}\partial_{2}h\int^{x_{2}}_{-r}\frac{t^{2}}{h^{3}(x_{1},t)}dt\\
&-\frac{3\omega_{1}}{2}\int_{\{|x'|<r\}}h\partial_{2}h\int^{x_{2}}_{-r}\frac{t^{2}}{h^{3}(x_{1},t)}dt.
\end{align*}
For simplicity of notations, we denote
\begin{align*}
\mathcal{K}_{1}=&\int_{\{|x'|<r\}}x_{2}\int^{x_{2}}_{-r}\frac{t^{2}}{h^{3}}dt,\quad\mathcal{K}_{2}=\int_{\{|x'|<r\}}\partial_{2}h\int^{x_{2}}_{-r}\frac{t^{2}}{h^{3}}dt,\\
\mathcal{K}_{3}=&\int_{\{|x'|<r\}}h\partial_{2}h\int^{x_{2}}_{-r}\frac{t^{2}}{h^{3}}dt.
\end{align*}
Similarly as (\ref{3.009}) and (\ref{3.0010}), a direct calculation yields that
\begin{align*}
\frac{\pi r^{2}}{8}\Phi_{33}^{(m)}(r/2;\varepsilon)\leq&\mathcal{K}_{1}+\frac{3\pi}{8}\Phi_{35}^{(m)}(r;\varepsilon)\leq\frac{\pi r^{2}}{2}\Phi_{33}^{(m)}(2r;\varepsilon),\\
\pi\left(\varepsilon+2^{-m}r^{m}\right)\Phi_{33}^{(m)}(r/2;\varepsilon)\leq&\mathcal{K}_{2}+\pi\Phi_{23}^{(m)}(r;\varepsilon)\leq\pi\left(\varepsilon+2^{\frac{m}{2}}r^{m}\right)\Phi_{33}^{(m)}(2r;\varepsilon),\\
\frac{\pi}{2}\left(\varepsilon+2^{-m}r^{m}\right)^{2}\Phi_{33}^{(m)}(r/2;\varepsilon)\leq&\mathcal{K}_{3}+\frac{\pi}{2}\Phi_{13}^{(m)}(r;\varepsilon)\leq\frac{\pi}{2}\left(\varepsilon+2^{\frac{m}{2}}r^{m}\right)^{2}\Phi_{33}^{(m)}(2r;\varepsilon).
\end{align*}

Thus, if $m=2$, we obtain
\begin{align}
\overline{T}^{(6)}_{1}\geq&-\frac{3\pi\mu\omega_{1}}{16}\left(6|\ln\varepsilon|+\frac{8Rr^{2}-r^{2}-2^{-4}r^{4}}{\varepsilon}\right),\label{21.020}\\
\overline{T}^{(6)}_{1}\leq&-\frac{3\pi\mu\omega_{1}}{16}\left(6|\ln\varepsilon|+\frac{Rr^{2}-4r^{2}-4r^{4}}{\varepsilon}\right);\label{21.021}
\end{align}
if $2<m<4$, then
\begin{align}
\overline{T}^{(6)}_{1}\geq&-\frac{3\pi\mu\omega_{1}}{8}\bigg(\frac{\Gamma_{34}^{(m)}}{2^{m-1}}\frac{2^{m+1}R+(2^{\frac{3}{2}m}-1)r^{m}}{\varepsilon^{2-{4}/{m}}}+\frac{3\Gamma_{36}^{(m)}}{\varepsilon^{3-\frac{6}{m}}}\notag\\
&\qquad\qquad\qquad+\frac{\Gamma_{34}^{(m)}}{2^{2m}}\frac{-2^{2m}r^{2}+2^{\frac{5}{2}m+2}Rr^{m}-r^{2m}}{\varepsilon^{3-{4}/{m}}}\bigg),\label{3.0014}\\
\overline{T}^{(6)}_{1}\leq&-\frac{3\pi\mu\omega_{1}}{8}\bigg(\frac{\Gamma_{34}^{(m)}}{2^{m-1}}\frac{2^{m+1}R-(2^{\frac{3}{2}m}-1)r^{m}}{\varepsilon^{2-{4}/{m}}}+\frac{3\Gamma_{36}^{(m)}}{\varepsilon^{3-\frac{6}{m}}}\notag\\
&\qquad\qquad\qquad+\frac{\Gamma_{34}^{(m)}}{2^{m-2}}\frac{-2^{m}r^{2}+Rr^{m}-2^{2m-2}r^{2m}}{\varepsilon^{3-{4}/{m}}}\bigg);\label{3.0015}
\end{align}
if $m=4$, then
\begin{align}
\overline{T}^{(6)}_{1}\geq&-\frac{3\pi\mu\omega_{1}}{8}\bigg(-2|\ln\varepsilon|+\frac{\Gamma_{34}^{(4)}}{8}\frac{32R+63r^{4}}{\varepsilon}\notag\\
&\qquad\qquad\qquad+\frac{3\Gamma_{36}^{(4)}}{\varepsilon^{\frac{3}{2}}}+\frac{\Gamma_{34}^{(4)}}{256}\frac{-256r^{2}+4096Rr^{4}-r^{8}}{\varepsilon^{2}}\bigg),\label{3.0016}\\
\overline{T}^{(6)}_{1}\leq&-\frac{3\pi\mu\omega_{1}}{8}\bigg(-2|\ln\varepsilon|+\frac{\Gamma_{34}^{(4)}}{8}\frac{32R-63r^{4}}{\varepsilon}\notag\\
&\qquad\qquad\qquad+\frac{3\Gamma_{36}^{(4)}}{\varepsilon^{\frac{3}{2}}}+\frac{\Gamma_{34}^{(4)}}{4}\frac{-16r^{2}+Rr^{4}-64r^{8}}{\varepsilon^{2}}\bigg);\label{3.0017}
\end{align}
if $m>4$, then
\begin{align}
\overline{T}^{(6)}_{1}\geq&-\frac{3\pi\mu\omega_{1}}{8}\bigg(\frac{2(\Gamma_{34}^{(m)}-2\Gamma_{24}^{(m)}-\frac{m+8}{m}\Gamma_{14}^{(m)})}{\varepsilon^{1-{4}/{m}}}+\frac{\Gamma_{34}^{(m)}}{2^{m-1}}\frac{2^{m+1}R+(2^{\frac{3}{2}m}-1)r^{m}}{\varepsilon^{2-{4}/{m}}}\notag\\
&\qquad\qquad\qquad+\frac{3\Gamma_{36}^{(m)}}{\varepsilon^{3-\frac{6}{m}}}+\frac{\Gamma_{34}^{(m)}}{2^{2m}}\frac{-2^{2m}r^{2}+2^{\frac{5}{2}m+2}Rr^{m}-r^{2m}}{\varepsilon^{3-{4}/{m}}}\bigg),\label{3.0018}\\
\overline{T}^{(6)}_{1}\leq&-\frac{3\pi\mu\omega_{1}}{8}\bigg(\frac{2(\Gamma_{34}^{(m)}-2\Gamma_{24}^{(m)}-\frac{m+8}{m}\Gamma_{14}^{(m)})}{\varepsilon^{1-{4}/{m}}}+\frac{\Gamma_{34}^{(m)}}{2^{m-1}}\frac{2^{m+1}R-(2^{\frac{3}{2}m}-1)r^{m}}{\varepsilon^{2-{4}/{m}}}\notag\\
&\qquad\qquad\qquad+\frac{3\Gamma_{36}^{(m)}}{\varepsilon^{3-\frac{6}{m}}}+\frac{\Gamma_{34}^{(m)}}{2^{m-2}}\frac{-2^{m}r^{2}+Rr^{m}-2^{2m-2}r^{2m}}{\varepsilon^{3-{4}/{m}}}\bigg);\label{3.0019}
\end{align}

Similarly, the above estimates (\ref{21.020})--(\ref{3.0019}) of $\overline{T}^{(6)}_{1}$ are still valid for $\overline{T}^{(6)}_{2}$ after replacing $\omega_{1}$ by $\omega_{2}$. In addition, by making use of the parity of integrand and the symmetry of domain, it follows that
\begin{align}\label{3.0020}
\overline{T}^{(6)}_{3}=0.
\end{align}

\subsection{Asymptotics of $\overline{\mathbf{F}}^{(0)}$ and $\overline{\mathbf{T}}^{(0)}$}

As before, we assume
\begin{align}\label{6.11}
\overline{\mathbf{u}}^{(0)}&=
\frac{1}{2}\mathbf{U}+\frac{1}{2}\boldsymbol{\omega}\times
\left(
\begin{array}{c}
          x_{1}    \\
          x_{2} \\
          \frac{|x'|^{m}}{2}-R
\end{array}
 \right).
\end{align}
Let $\overline{p}^{(0)}=C$. It is straightforward to check
\begin{align}\label{6.12}
\overline{\mathbf{F}}^{(0)}&=\mathbf{0}+O(\mathbf{1}),\quad\overline{\mathbf{T}}^{(0)}=\mathbf{0}+O(\mathbf{1}),
\end{align}
which implies this type of the fluid flow contributes only to the $O(\mathbf{1})$ order of the forces.

\subsection{Justification}
Utilizing integration by parts for $I_{\Omega}[\mathbf{u}]$, in which $\mathbf{u}$ is the solution of $(\ref{11.3})$, we obtain
$$I_{\Omega}[\mathbf{u}]
=-\frac{\mathbf{U}\cdot\mathbf{F}+\boldsymbol{\omega}\cdot\mathbf{T}}{2},$$
which indicates the minimal value $I_{\Omega}[\mathbf{u}]$ can be expressed as the terms of $\mathbf{F}$ and $\mathbf{T}$.

Denote by
\begin{align*}
\mathbf{A}:\mathbf{B}=\mathrm{tr}(\mathbf{A}\mathbf{B})=\sum\limits^{3}_{i,j=1}\mathbf{A}_{ij}\mathbf{B}_{ij}
\end{align*}
the trace product of two symmetric tensors $\mathbf{A}$ and $\mathbf{B}$. Notice that the fluid stress $\boldsymbol{\sigma}(\mathbf{u})$ corresponding to $\mathbf{u}$ gives a maximum for the following functional:
\begin{align}\label{11.78}
\boldsymbol{\sigma}(\mathbf{u})&=\mathrm{arg}\max_{\mathbf{S}\in\mathcal{A}^{\ast}}I^{\ast}_{\Omega}[\mathbf{S}],
\end{align}
where
\begin{align*}
I^{\ast}_{\Omega}[\mathbf{S}]&=\int_{\partial D_{1}}\mathbf{u}\cdot\mathbf{S}\mathbf{n}\,dS-\frac{1}{4\mu}\int_{\Omega}\left(\mathrm{tr}\mathbf{S}^{2}-\frac{(\mathrm{tr}\mathbf{S})^{2}}{3}\right)dx,\\
\mathcal{A}^{\ast}&=\big\{\mathbf{S}\in\mathbb{R}^{3\times3}:~S_{ij}\in L^{2}(\Omega),~\mathbf{S}=\mathbf{S}^{T},~\nabla\cdot\mathbf{S}=\mathbf{0}~\mathrm{in}~\Omega,~\mathbf{S}\mathbf{n}=\mathbf{0}~\mathrm{on}~\partial\Omega\big\}.
\end{align*}
We call problem $(\ref{11.78})$ the dual variational principle corresponding to $(\ref{5.3})$. To obtain the justification in Theorem $\ref{thm1}$, we need to find a test vector $\boldsymbol{\mathfrak{U}}\in\mathcal{A}$ and a test tensor $\boldsymbol{\mathfrak{S}}\in\mathcal{A}^{\ast}$ satisfying
\begin{align*}
\mathrm{Err}:=&\left|I_{\Omega_{r}}[\boldsymbol{\mathfrak{U}}]-I^{\ast}_{\Omega_{r}}[\boldsymbol{\mathfrak{S}}]\right|\notag\\
=&\left|\frac{1}{2}\int_{\Omega}\boldsymbol{\sigma}(\boldsymbol{\mathfrak{U}}):\mathbf{D}(\boldsymbol{\mathfrak{U}})dx-\int_{\partial D_{1}}\boldsymbol{\mathfrak{U}}\cdot\boldsymbol{\mathfrak{S}}\mathbf{n}~dS+\frac{1}{4\mu}\int_{\Omega}\left(\mathrm{tr}\boldsymbol{\mathfrak{S}}^{2}-\frac{(\mathrm{tr}\boldsymbol{\mathfrak{S}})^{2}}{3}\right)dx\right|\notag\\
=&O(1).
\end{align*}
$\mathrm{\mathbf{Step}}~\mathrm{\mathbf{1}}$. The construction of test vector $\boldsymbol{\mathfrak{U}}\in\mathcal{A}$.

Choose
\begin{align}\label{11.81}
\boldsymbol{\mathfrak{U}}=
\begin{cases}
\sum\limits^{6}_{i=0}\overline{\mathbf{u}}^{(i)},&~~~~~\Omega_{r},\\
\widetilde{\mathbf{u}},&~~~~~\Omega\setminus\Omega_{r},
\end{cases}
\end{align}
where $\overline{\mathbf{u}}^{(i)}$, $i\in\{0,\cdots,6\}$, are defined by (\ref{6.11}), (\ref{11.18})--(\ref{11.19}), (\ref{11.29})--(\ref{11.30}), (\ref{11.33})--(\ref{11.344}), (\ref{11.52})--(\ref{11.533}), (\ref{11.55})--(\ref{15.02}) and (\ref{11.60})--(\ref{11.61}). By virtue of the Kirszbraun theorem \cite{F1969}, we extend $\sum\limits^{6}_{i=0}\overline{\mathbf{u}}^{(i)}$ of $\Omega_{r}$ to the part of the domain $\Omega\setminus\Omega_{r}$ to obtain the function $\widetilde{\mathbf{u}}$ in $(\ref{11.81})$, which satisfies $I_{\Omega\setminus\Omega_{r}}[\widetilde{\mathbf{u}}]=O(1)$. Then
\begin{equation}\label{11.82}
I_{\Omega}[\boldsymbol{\mathfrak{U}}]=I_{\Omega_{r}}[\boldsymbol{\mathfrak{U}}]+O(1).
\end{equation}
$\mathrm{\mathbf{Step}}~\mathrm{\mathbf{2}}$. The construction of test tensor $\boldsymbol{\mathfrak{S}}\in\mathcal{A}^{\ast}$.

Let
\begin{align}\label{11.83}
\boldsymbol{\mathfrak{S}}=
\begin{cases}
\sum\limits^{6}_{i=0}\mathbf{S}^{(i)},&~~~~~\Omega_{r},\\
\mathbf{0},&~~~~~\Omega\setminus\Omega_{r},
\end{cases}
\end{align}
where $\mathbf{S}^{(i)}$, $i\in\{0,\cdots,6\}$, correspond to the dual variational formulation to $(\ref{5.1})$, defined by
\begin{align}\label{11.84}
\boldsymbol{\sigma}(\mathbf{u}^{(i)})&=\mathrm{arg}\max_{\mathbf{S}\in\mathcal{A}^{\ast}_{\Omega_{r}}}I^{\ast}_{\Omega_{r}}[\mathbf{S}]
\end{align}
with the functional $I^{\ast}_{\Omega_{r}}[\cdot]$ having the similar form as the one in $(\ref{11.78})$ and $$\mathcal{A}^{\ast}_{\Omega_{r}}=\big\{\mathbf{S}\in\mathbb{R}^{3\times3}:~S_{ij}\in L^{2}(\Omega_{r}),~\mathbf{S}=\mathbf{S}^{T},~\nabla\cdot\mathbf{S}=\mathbf{0}~\mathrm{in}~\Omega_{r},~\mathbf{S}\mathbf{n}=\mathbf{0}~\mathrm{on}~\partial\Omega_{r}\big\}.$$
Combining with $(\ref{11.82})$, it follows that
\begin{align}\label{11.85}
\mathrm{Err}&\leq\left|I_{\Omega_{r}}[\boldsymbol{\mathfrak{U}}]-I^{\ast}_{\Omega_{r}}[\boldsymbol{\mathfrak{S}}]\right|+O(1).
\end{align}

Write
\begin{align*}
\Omega_{\frac{r}{2}}:=\left\{(x',x_{3})\in\mathbb{R}^{3}:\,-h/2<x_{3}<h/2,~|x'|<r/4\right\},
\end{align*}
and let $\mathbf{S}^{(i)}=\mathbf{0}$ in $\Omega_{r}\setminus\Omega_{\frac{r}{2}}$ so that $\mathbf{S}^{(i)}$ satisfies the traction-free condition on $\partial\Omega_{r}$. In order to construct $\mathbf{S}^{(i)}$ satisfying
\begin{align}\label{5.5}
\nabla\cdot\mathbf{S}^{(i)}=\mathbf{0},\quad\mathrm{in}\;\Omega_{\frac{r}{2}},
\end{align}
we retain the terms of $\boldsymbol{\sigma}(\mathbf{u}^{(i)})$ contributing to the leading terms of asymptotics, and add/subtract correcting terms to the entries.

Thus we correct $(\ref{11.83})$ by
\begin{align}\label{11.86}
\boldsymbol{\mathfrak{S}}=
\begin{cases}
\sum\limits^{6}_{i=0}\mathbf{S}^{(i)},&~~~~~\Omega_{\frac{r}{2}},\\
\mathbf{0},&~~~~~\Omega\setminus\Omega_{\frac{r}{2}},
\end{cases}
\end{align}
which leads to an adjustment of the error $(\ref{11.85})$ as follows:
\begin{align}\label{11.87}
\mathrm{Err}\leq\left|I_{\Omega_{\frac{r}{2}}}[\boldsymbol{\mathfrak{U}}]-I^{\ast}_{\Omega_{\frac{r}{2}}}[\boldsymbol{\mathfrak{S}}]\right|+O(1).
\end{align}

Observe that
\begin{equation}\label{11.88}
I_{\Omega_{\frac{r}{2}}}[\boldsymbol{\mathfrak{U}}]-I^{\ast}_{\Omega_{\frac{r}{2}}}[\boldsymbol{\mathfrak{S}}]=\mu\int_{\Omega_{\frac{r}{2}}}\mathrm{tr}\left[\mathbf{D}(\boldsymbol{\mathfrak{U}})-\frac{1}{2\mu}\left(\boldsymbol{\mathfrak{S}}-\frac{\mathrm{tr}\boldsymbol{\mathfrak{S}}}{3}\mathbf{E}\right)\right]^{2}dx,
\end{equation}
where $\mathbf{E}\in\mathbb{R}^{3\times3}$ denotes the unit tensor. According to $(\ref{11.81})$ and $(\ref{11.86})$, error $(\ref{11.88})$ can be expressed as
\begin{align}\label{11.89}
I_{\Omega_{\frac{r}{2}}}[\boldsymbol{\mathfrak{U}}]-I^{\ast}_{\Omega_{\frac{r}{2}}}[\boldsymbol{\mathfrak{S}}]=\sum\limits^{6}_{i=0}\ell[i,i]+2\sum_{i<j}\ell[i,j],
\end{align}
where $$\ell[i,j]=\mu\int_{\Omega_{\frac{r}{2}}}\left[\mathbf{D}(\overline{\mathbf{u}}^{(i)})-\frac{1}{2\mu}\left(\mathbf{S}^{(i)}-\frac{\mathrm{tr}\mathbf{S}^{(i)}}{3}\mathbf{E}\right)\right]:\left[\mathbf{D}(\overline{\mathbf{u}}^{(j)})-\frac{1}{2\mu}\left(\mathbf{S}^{(j)}-\frac{\mathrm{tr}\mathbf{S}^{(j)}}{3}\mathbf{E}\right)\right].$$
Next we select the test tensors $\mathbf{S}^{(i)}$ satisfying $\ell[i,j]=O(1)$ for all $i,j=0,...,6$.
$\mathrm{\mathbf{Step}}~\mathrm{\mathbf{3}}$. Selection of particular $\mathbf{S}^{(i)}\in\mathcal{A}^{\ast}_{\Omega_{r}}$.

Since $\overline{\mathbf{u}}^{(0)}$, $\overline{\mathbf{u}}^{(4)}$ and $\overline{\mathbf{u}}^{(5)}$ contribute only to the $O(\mathbf{1})$ order in the asymptotics of $\overline{\mathbf{F}}$ and $\overline{\mathbf{T}}$, we choose
\begin{equation}\label{11.90}
\mathbf{S}^{(0)}=\mathbf{S}^{(4)}=\mathbf{S}^{(5)}=\mathbf{0},\quad\mathrm{in}\;\Omega_{r}.
\end{equation}
It is clear that
\begin{align}\label{11.91}
\ell[0,0]=O(1),\quad\ell[4,4]=O(1),\quad\ell[5,5]=O(1).
\end{align}
When $i=1$, we let
\begin{align}\label{11.92}
\mathbf{S}^{(1)}&=\mu(U_{1}-\omega_{2}R)
\left(
\begin{array}{ccc}
          0 & 0 & H(x') \\
          0 & 0 &  0     \\
          H(x') & 0 & -B(x')x_{3}
\end{array}
 \right),\quad\mathrm{in}~\Omega_{\frac{r}{2}},
\end{align}
where $B$ and $H$ are defined by $(\ref{11.19})$. Set $\mathbf{S}^{(1)}=\mathbf{0}~\mathrm{in}~\Omega_{r}\setminus\Omega_{\frac{r}{2}}$. A direct calculation yields
\begin{align}\label{11.93}
\ell[1,1]=O(1).
\end{align}
Similarly as before, for $i=2$, we select
\begin{align}\label{11.95}
\mathbf{S}^{(2)}&=\mu(U_{2}+\omega_{1}R)
\left(
\begin{array}{ccc}
          0 & 0 &  0      \\
          0 & 0 & H(x') \\
          0 & H(x') &-B(x')x_{3}
\end{array}
 \right),\quad\mathrm{in}\;\Omega_{\frac{r}{2}},
\end{align}
with $B$ and $H$ defined by $(\ref{11.30})$. Pick $\mathbf{S}^{(2)}=\mathbf{0}~\mathrm{in}~\Omega_{r}\setminus\Omega_{\frac{r}{2}}$. Similarly, we have
\begin{equation}\label{11.96}
\ell[2,2]=O(1).
\end{equation}

Since there are many terms contributing to the leading orders in the corresponding $\boldsymbol{\sigma}(\overline{\mathbf{u}}^{(i)}),\,i\in\{3,6\}$, we construct $\mathbf{S}^{(i)}$ by correcting $\boldsymbol{\sigma}(\overline{\mathbf{u}}^{(i)})$ according to (\ref{5.5}). Then, for $i\in\{3,6\}$, we set
\begin{align}\label{11.97}
\mathbf{S}_{kl}^{(i)}=
\begin{cases}
2\mu\partial_{k}\overline{\mathbf{u}}_{k}^{(i)}-\overline{p}^{(i)}+q^{(i)}_{k},&~~~~~k=l\\
\mu(\partial_{l}\overline{\mathbf{u}}^{(i)}_{k}+\partial_{k}\overline{\mathbf{u}}_{l}^{(i)}),&~~~~~k\neq l
\end{cases}\quad k,l\in\{1,2,3\},
\end{align}
where $\overline{\mathbf{u}}^{(3)}$ and $\overline{p}^{(3)}$ are determined by (\ref{11.33})--(\ref{11.344}), $\overline{\mathbf{u}}^{(6)}$ and $\overline{p}^{(6)}$ are defined by (\ref{11.60})--(\ref{11.62}), and $q_{k}^{(i)},\,k\in\{1,2,3\}$ are given by
\begin{align}
q_{1}^{(i)}&=\int^{x_{1}}_{0}\alpha^{(i)}\left[\Delta A_{1}-\partial_{1}A_{3}+3x_{3}^{2}(\Delta B_{1}+\partial_{1}B_{3})\right]dx_{1},\label{11.98}\\
q_{2}^{(i)}&=\int^{x_{2}}_{0}\alpha^{(i)}\left[\Delta A_{2}-\partial_{2}A_{3}+3x_{3}^{2}(\Delta B_{2}+\partial_{2}B_{3})\right]dx_{2},\label{11.988}\\
q_{3}^{(i)}&=-\alpha^{(i)}\left(\frac{\Delta A_{3}}{2}x_{3}^{2}+\frac{\Delta B_{3}}{4}x_{3}^{4}\right),\label{11.9888}\\
\alpha^{(i)}&=
\begin{cases}
\mu U_{3},\quad\quad\quad i=3,\\
\mu,\quad\quad\quad\quad i=6,
\end{cases}\notag
\end{align}
where $A_{j},B_{j},\,j\in\{1,2,3\}$ are defined by (\ref{11.34})--(\ref{11.344}) for $i=3$ and (\ref{11.61}) for $i=6$, respectively, and $\Delta:=\partial_{x_{1}x_{1}}+\partial_{x_{2}x_{2}}$. Similarly, take $\mathbf{S}^{(i)}=\mathbf{0}~\mathrm{in}~\Omega_{r}\setminus\Omega_{\frac{r}{2}},\,i\in\{3,6\}$. By evaluating the corresponding functional difference, we arrive at
\begin{align}\label{11.99}
\ell[i,i]=O(1),\quad i\in\{3,6\}.
\end{align}

Finally, it follows from $(\ref{11.90})$, $(\ref{11.92})$, $(\ref{11.95})$ and (\ref{11.97})--(\ref{11.9888}) that, for $i,j=0,\cdots,6$,
\begin{align}\label{11.100}
\ell[i,j]&=O(1),\quad i<j.
\end{align}
In light of $(\ref{11.89})$, $(\ref{11.91})$, $(\ref{11.93})$, $(\ref{11.96})$, $(\ref{11.99})$ and $(\ref{11.100})$, we obtain
$$\mathrm{Err}=O(1).$$

Theorem $\ref{thm1}$ is therefore proved by using (\ref{21.010})--(\ref{11.28}), (\ref{21.012})--(\ref{11.32}), (\ref{11.51}), (\ref{11.54}), (\ref{21.015}), (\ref{21.017})--(\ref{3.0020}) and (\ref{6.12}).

\section{Proof of Theorem \ref{thm3}}
For simplicity, we only give the calculations of $\overline{\mathbf{F}}^{(i)}$ and $\overline{\mathbf{T}}^{(i)}$, $i=3,6$. Other parts are similar and thus omitted. First, observe that under the hypothesis (\ref{2.001}), the unit normal $\mathbf{n}$ to $\Sigma^{+}_{r}$ is given by
\begin{align}\label{665.01}
\mathbf{n}&=
\begin{cases}
(0,0,-1),&\Sigma^{+}_{s},\\
\frac{1}{\sqrt{1+(|x'|-s)^{2}}}\left(\dfrac{x_{1}(|x'|-s)}{|x'|},\dfrac{x_{2}(|x'|-s)}{|x'|},-1\right),&\Sigma^{+}_{r}\setminus\Sigma^{+}_{s},
\end{cases}
\end{align}
and the corresponding $\boldsymbol{\nu}$ is given by
\begin{align}\label{665.02}
\boldsymbol{\nu}&=
\begin{cases}
(x_{1},x_{2},-R),&\Sigma^{+}_{s},\\
\Big(x_{1},x_{2},\frac{(|x'|-s)^{2}}{2}-R\Big),&\Sigma^{+}_{r}\setminus\Sigma^{+}_{s}.
\end{cases}
\end{align}
Throughout this section, for $s,r>0$, $i>0$ and $j\geq0$, denote
\begin{align*}
\Psi_{ij}(r;\varepsilon)=&\int_{0}^{r}\frac{(t+s)^{j}}{(\varepsilon+t^{2})^{i}}\,dt.
\end{align*}
\subsection{Asymptotics of $\overline{\mathbf{F}}^{(3)}$ and $\overline{\mathbf{T}}^{(3)}$}
In light of (\ref{11.39})--(\ref{11.42}), we have
\begin{align}
&\int_{\Sigma^{+}_{r}}\left((\partial_{1}A_{3})x_{3}+(\partial_{1}B_{3})x_{3}^{3}-6B_{1}x_{3}\right)n_{1}\,dS=3\pi\Psi_{11}(r-s;\varepsilon),\label{111.43}\\
&\int_{\Sigma^{+}_{r}}\left((\partial_{2}A_{3})x_{3}+(\partial_{2}B_{3})x_{3}^{3}-6B_{2}x_{3}\right)n_{2}\,dS=3\pi\Psi_{11}(r-s;\varepsilon),\label{111.44}\\
&\int_{\Sigma^{+}_{r}}2(A_{3}+3B_{3}x_{3}^{2})n_{3}\,dS=-12\pi\Psi_{11}(r-s;\varepsilon),\label{111.45}
\end{align}
and
$$\int_{\Sigma^{+}_{r}}-\frac{\overline{p}^{(3)}}{\mu U_{3}}n_{3}\,dS=\int_{\Sigma^{+}_{s}}-\frac{\overline{p}^{(3)}}{\mu U_{3}}n_{3}\,dS+\int_{\Sigma^{+}_{r}\setminus\Sigma^{+}_{s}}-\frac{\overline{p}^{(3)}}{\mu U_{3}}n_{3}\,dS
:=K_{1}+K_{2},$$
where
\begin{align}
K_{1}=&-\frac{3\pi s^{2}}{\varepsilon}-\frac{3\pi(2s^{2}r^{2}-s^{4})}{2\varepsilon^{3}},\label{111.50}\\
K_{2}=&6\pi\Psi_{11}(r-s;\varepsilon)+6\pi s^{2}\Psi_{31}(r-s;\varepsilon)-6\pi\Psi_{33}(r-s;\varepsilon).\label{111.46}
\end{align}

In light of \eqref{GBD03} and combining with (\ref{111.43})--(\ref{111.46}), we obtain
\begin{align}
\overline{\mathbf{F}}^{(3)}&=-\frac{3\pi\mu U_{3}}{2}\left(\frac{1}{\varepsilon}+\frac{6s\Gamma_{21}^{(2)}}{\varepsilon^{\frac{3}{2}}}+\frac{8s^{2}\Gamma_{32}^{(2)}}{\varepsilon^{2}}+\frac{2r^{2}s^{2}-s^{4}}{\varepsilon^{3}}\right)\mathbf{e}_{3}+O(\mathbf{1}),\label{111.51}\\
\overline{\mathbf{T}}^{(3)}&=\mathbf{0}+O(\mathbf{1}).\label{111.52}
\end{align}

\subsection{Estimates of $\overline{\mathbf{F}}^{(6)}$ and $\overline{\mathbf{T}}^{(6)}$}
According to (\ref{11.60})--(\ref{3.008}), we now estimate each component of $\overline{\mathbf{F}}^{(6)}$ and $\overline{\mathbf{T}}^{(6)}$, respectively.  As above, $O(1)$ is also omitted in the following estimates of $\overline{F}^{(6)}_{i}$ and $\overline{T}^{(6)}_{i}$ for $i=1,2,3$.

$\mathrm{\mathbf{Step}}~\mathrm{\mathbf{1}}$. Estimates of $\overline{F}^{(6)}_{i},~i=1,2,3$.

On one hand, by integrating by parts, we obtain
\begin{align}
\int_{\Sigma^{+}_{r}}-\frac{\overline{p}^{(6)}}{3\mu}n_{1}=&\int_{\{|x'|<r\}}\partial_{1}hG_{1}\notag\\
\leq&\omega_{2}\int^{r}_{-r}dx_{2}\int^{r}_{-r}\partial_{1}h\int^{x_{1}}_{r}\frac{t^{2}}{h^{3}(t,x_{2})}dt~dx_{1}\notag\\
\leq&\omega_{2}\left(\varepsilon+(\sqrt{2}r-s)^{2}\right)\int_{\{|x'|<2r\}}\frac{x_{1}^{2}}{h^{3}}-\omega_{2}\int_{\{|x'|<r\}}\frac{x_{1}^{2}}{h^{2}}\notag\\
\leq&\omega_{2}\left(\varepsilon+(\sqrt{2}r-s)^{2}\right)\left(\pi\Psi_{33}(2r-s;\varepsilon)+\frac{\pi s^{4}}{4\varepsilon^{3}}\right)\notag\\
&-\omega_{2}\left(\pi\Psi_{23}(r-s;\varepsilon)+\frac{\pi s^{4}}{4\varepsilon^{2}}\right).\label{33.009}
\end{align}
On the other hand, applying integration by parts and making use of $0<s<(\sqrt{2}-1)r$, we have
\begin{align}
\int_{\Sigma^{+}_{r}}-\frac{\overline{p}^{(6)}}{3\mu}n_{1}=&\int_{\{|x'|<r\}}\partial_{1}hG_{1}\notag\\
\geq&\omega_{2}\int^{\frac{r+s}{2}}_{-\frac{r+s}{2}}dx_{2}\int^{\frac{r+s}{2}}_{-\frac{r+s}{2}}\partial_{1}h\int^{x_{1}}_{r}\frac{t^{2}}{h^{3}(t,x_{2})}dt~dx_{1}\notag\\
\geq&\omega_{2}\left(\varepsilon+2^{-2}(r-s)^{2}\right)\int_{\{|x'|<\frac{r+s}{2}\}}\frac{x_{1}^{2}}{h^{3}}-\omega_{2}\int_{\{|x'|<r\}}\frac{x_{1}^{2}}{h^{2}}\notag\\
\geq&\omega_{2}\left(\varepsilon+2^{-2}(r-s)^{2}\right)\left(\pi\Psi_{33}\big((r-s)/2;\varepsilon\big)+\frac{\pi s^{4}}{4\varepsilon^{3}}\right)\notag\\
&-\omega_{2}\left(\pi\Psi_{23}(r-s;\varepsilon)+\frac{\pi s^{4}}{4\varepsilon^{2}}\right).\label{33.0010}
\end{align}
In addition, by using (\ref{11.644}), we obtain
\begin{align}
\int_{\Sigma^{+}_{r}}\left((\partial_{1}A_{3})x_{3}+(\partial_{1}B_{3})x_{3}^{3}-6B_{1}x_{3}\right)n_{3}
&=\frac{3\pi\omega_{2}s^{4}}{4\varepsilon^{2}}+3\pi\omega_{2}\Psi_{23}(r-s;\varepsilon).\label{111.67}
\end{align}
Thus, in view of \eqref{GBD06} and utilizing (\ref{33.009})--(\ref{111.67}), we arrive at
\begin{align*}
\overline{F}^{(6)}_{1}\geq&\frac{3\pi\mu\omega_{2}}{16}\bigg(\frac{24s\Gamma^{(2)}_{21}}{\varepsilon^{\frac{1}{2}}}+\frac{24s^{2}\Gamma_{32}^{(2)}+(r-s)^{2}}{\varepsilon}+\frac{8s^{3}(\Gamma_{31}^{(2)}+\Gamma_{21}^{(2)})+6s(r-s)^{2}\Gamma_{21}^{(2)}}{\varepsilon^{\frac{3}{2}}}\notag\\
&+\frac{6s^{2}(r-s)^{2}\Gamma_{32}^{(2)}+4s^{4}}{\varepsilon^{2}}+\frac{2s^{3}(r-s)^{2}\Gamma_{31}^{(2)}}{\varepsilon^{\frac{5}{2}}}+\frac{s^{4}(r-s)^{2}}{\varepsilon^{3}}\bigg),\\
\overline{F}^{(6)}_{1}\leq&\frac{3\pi\mu\omega_{2}}{4}\bigg(\frac{6s\Gamma_{21}^{(2)}}{\varepsilon^{\frac{1}{2}}}+\frac{6s^{2}\Gamma_{32}^{(2)}+(\sqrt{2}r-s)^{2}}{\varepsilon}+\frac{2s^{3}(\Gamma_{31}^{(2)}+\Gamma_{21}^{(2)})+6s(\sqrt{2}r-s)^{2}\Gamma_{21}^{(2)}}{\varepsilon^{\frac{3}{2}}}\notag\\
&+\frac{6s^{2}(\sqrt{2}r-s)^{2}\Gamma_{32}^{(2)}+s^{4}}{\varepsilon^{2}}+\frac{2s^{3}(\sqrt{2}r-s)^{2}\Gamma_{31}^{(2)}}{\varepsilon^{\frac{5}{2}}}+\frac{s^{4}(\sqrt{2}r-s)^{2}}{\varepsilon^{3}}\bigg).
\end{align*}

Similarly as above, in view of (\ref{11.65}) and (\ref{3.008}), we obtain
\begin{align*}
\overline{F}^{(6)}_{2}\geq&-\frac{3\pi\mu\omega_{1}}{4}\bigg(\frac{6s\Gamma_{21}^{(2)}}{\varepsilon^{\frac{1}{2}}}+\frac{6s^{2}\Gamma_{32}^{(2)}+(\sqrt{2}r-s)^{2}}{\varepsilon}+\frac{2s^{3}(\Gamma_{31}^{(2)}+\Gamma_{21}^{(2)})+6s(\sqrt{2}r-s)^{2}\Gamma_{21}^{(2)}}{\varepsilon^{\frac{3}{2}}}\notag\\
&+\frac{6s^{2}(\sqrt{2}r-s)^{2}\Gamma_{32}^{(2)}+s^{4}}{\varepsilon^{2}}+\frac{2s^{3}(\sqrt{2}r-s)^{2}\Gamma_{31}^{(2)}}{\varepsilon^{\frac{5}{2}}}+\frac{s^{4}(\sqrt{2}r-s)^{2}}{\varepsilon^{3}}\bigg),\\
\overline{F}^{(6)}_{2}\leq&-\frac{3\pi\mu\omega_{1}}{16}\bigg(\frac{24s\Gamma^{(2)}_{21}}{\varepsilon^{\frac{1}{2}}}+\frac{24s^{2}\Gamma_{32}^{(2)}+(r-s)^{2}}{\varepsilon}+\frac{8s^{3}(\Gamma_{31}^{(2)}+\Gamma_{21}^{(2)})+6s(r-s)^{2}\Gamma_{21}^{(2)}}{\varepsilon^{\frac{3}{2}}}\notag\\
&+\frac{6s^{2}(r-s)^{2}\Gamma_{32}^{(2)}+4s^{4}}{\varepsilon^{2}}+\frac{2s^{3}(r-s)^{2}\Gamma_{31}^{(2)}}{\varepsilon^{\frac{5}{2}}}+\frac{s^{4}(r-s)^{2}}{\varepsilon^{3}}\bigg),\\
\overline{F}^{(6)}_{3}\geq&\frac{3\pi\mu(r+s)(\omega_{1}+\omega_{2})}{2}\bigg(\frac{1}{2\varepsilon}+\frac{3s\Gamma_{21}^{(2)}}{\varepsilon^{\frac{3}{2}}}+\frac{3s^{2}\Gamma_{32}^{(2)}}{\varepsilon^{2}}+\frac{s^{3}\Gamma_{31}^{(2)}}{\varepsilon^{\frac{5}{2}}}+\frac{s^{4}}{2\varepsilon^{3}}\bigg),\\
\overline{F}^{(6)}_{3}\leq&3\pi\mu r(\omega_{1}+\omega_{2})\bigg(\frac{1}{2\varepsilon}+\frac{3s\Gamma_{21}^{(2)}}{\varepsilon^{\frac{3}{2}}}+\frac{3s^{2}\Gamma_{32}^{(2)}}{\varepsilon^{2}}+\frac{s^{3}\Gamma_{31}^{(2)}}{\varepsilon^{\frac{5}{2}}}+\frac{s^{4}}{2\varepsilon^{3}}\bigg).
\end{align*}

$\mathrm{\mathbf{Step}}~\mathrm{\mathbf{2}}$. Estimates of $\overline{T}^{(6)}_{i},~i=1,2,3$.

First, combining with (\ref{11.644})--(\ref{3.008}) and (\ref{665.01})--(\ref{665.02}), we have
\begin{align*}
&\int_{\Sigma^{+}_{r}}\left((\partial_{2}A_{3})x_{3}+(\partial_{2}B_{3})x_{3}^{3}-6B_{2}x_{3}\right)n_{3}\nu_{3}\,dS\\
=&3\pi\omega_{1}R\Psi_{33}(r-s;\varepsilon)-\frac{9\pi\omega_{1}s^{2}}{2}\Phi_{23}^{(2)}(r-s;\varepsilon)-\frac{3\pi\omega_{1}s^{3}}{2}\Phi_{22}^{(2)}(r-s;\varepsilon)+\frac{3\pi\omega_{1}Rs^{4}}{4\varepsilon^{2}},
\end{align*}
and
\begin{align*}
\int_{\Sigma^{+}_{r}}\frac{\overline{p}^{(6)}}{\mu}n_{3}\nu_{2}=&3\omega_{1}\left(\pi\Psi_{13}(r-s;\varepsilon)+\frac{\pi s^{4}}{4\varepsilon}\right)-6\omega_{1}\int_{\{|x'|<r\}}x_{2}\int^{x_{2}}_{-r}\frac{t^{2}}{h^{3}}dt,\\
\int_{\Sigma^{+}_{r}}-\frac{\overline{p}^{(6)}}{\mu}n_{2}\nu_{3}=&\frac{6\omega_{1}R+3\omega_{1}\varepsilon}{2}\int_{\{|x'|<r\}}\partial_{2}h\int^{x_{2}}_{-r}\frac{t^{2}}{h^{3}(x_{1},t)}dt\\
&-\frac{3\omega_{1}}{2}\int_{\{|x'|<r\}}h\partial_{2}h\int^{x_{2}}_{-r}\frac{t^{2}}{h^{3}(x_{1},t)}dt.
\end{align*}
Similarly as before, we denote
\begin{align*}
\mathcal{K}_{1}=&\int_{\{|x'|<r\}}x_{2}\int^{x_{2}}_{-r}\frac{t^{2}}{h^{3}}dt,\quad\mathcal{K}_{2}=\int_{\{|x'|<r\}}\partial_{2}h\int^{x_{2}}_{-r}\frac{t^{2}}{h^{3}}dt,\\
\mathcal{K}_{3}=&\int_{\{|x'|<r\}}h\partial_{2}h\int^{x_{2}}_{-r}\frac{t^{2}}{h^{3}}dt.
\end{align*}
Then it follows from integration by parts and $0<s<(\sqrt{2}-1)r$ that
\begin{align*}
\mathcal{K}_{1}\leq&\frac{\pi r^{2}}{2}\Psi_{33}(2r-s;\varepsilon)-\frac{3\pi}{8}\Psi_{35}(r-s;\varepsilon)+\frac{2\pi r^{2}s^{4}-\pi s^{6}}{16\varepsilon^{3}},\\
\mathcal{K}_{1}\geq&\frac{\pi(r+s)^{2}}{8}\Psi_{33}((r-s)/2;\varepsilon)-\frac{3\pi}{8}\Psi_{35}(r-s;\varepsilon)+\frac{\pi(r+s)^{2}s^{4}-2\pi s^{6}}{32\varepsilon^{3}},\\
\mathcal{K}_{2}\leq&\left(\varepsilon+(\sqrt{2}r-s)^{2}\right)\left(\pi\Psi_{33}(2r-s;\varepsilon)+\frac{\pi s^{4}}{4\varepsilon^{3}}\right)-\left(\pi\Psi_{23}(r-s;\varepsilon)+\frac{\pi s^{4}}{4\varepsilon^{2}}\right),\\
\mathcal{K}_{2}\geq&\left(\varepsilon+2^{-2}(r-s)^{2}\right)\left(\pi\Psi_{33}((r-s)/2;\varepsilon)+\frac{\pi s^{4}}{4\varepsilon^{3}}\right)-\left(\pi\Psi_{23}(r-s;\varepsilon)+\frac{\pi s^{4}}{4\varepsilon^{2}}\right),\\
\mathcal{K}_{3}\leq&\frac{\pi\left(\varepsilon+(\sqrt{2}r-s)^{2}\right)^{2}}{2}\left(\Psi_{33}(2r-s;\varepsilon)+\frac{\pi s^{4}}{4\varepsilon^{3}}\right)-\frac{\pi s^{4}}{8\varepsilon},\\
\mathcal{K}_{3}\geq&\frac{\pi\left(\varepsilon+2^{-2}(r-s)^{2}\right)^{2}}{2}\left(\Psi_{33}((r-s)/2;\varepsilon)+\frac{\pi s^{4}}{4\varepsilon^{3}}\right)-\frac{\pi s^{4}}{8\varepsilon}.
\end{align*}
Therefore, a direct calculation yields
\begin{align*}
\overline{T}^{(6)}_{1}\geq&-\frac{3\pi\mu\omega_{1}}{32}\bigg(\frac{4R(\sqrt{2}r-s)^{2}-(r+s)^{2}-2^{-4}(r-s)^{4}}{\varepsilon}\\
&+\frac{2s^{4}\big(4R(\sqrt{2}r-s)^{2}+(s-r)^{2}-2r^{2}-2^{-4}(r-s)^{4}\big)}{\varepsilon^{3}}\bigg),\\
\overline{T}^{(6)}_{1}\leq&-\frac{3\pi\mu\omega_{1}}{16}\bigg(\frac{2\big(R(r-s)^{2}-4r^{2}-(\sqrt{2}r-s)^{4}\big)}{\varepsilon}\\
&+\frac{s^{4}\big(R(r-s)^{2}-4r^{2}+2s^{2}-(\sqrt{2}r-s)^{4}\big)}{\varepsilon^{3}}\bigg).
\end{align*}

Similarly, the corresponding estimates of $\overline{T}^{(6)}_{2}$ are the same as $\overline{T}^{(6)}_{1}$ except replacing $\omega_{1}$ by $\omega_{2}$. Moreover, it is easy to verify that
\begin{align*}
\overline{T}^{(6)}_{3}=0.
\end{align*}

\section{Approximation method in $2\mathrm{D}$}
The variational formula in $2\mathrm{D}$ is similar to the case of $3\mathrm{D}$ in section 3.1 and thus omitted here.
\subsection{Problem in the small gap width between $D_{1}$ and $D_{2}$ in $2\mathrm{D}$}

Similarly as $3\mathrm{D}$, denote by
\begin{equation*}
h(x_{1})=\varepsilon+|x_{1}|^{m}
\end{equation*}
the vertical distance between particles $D_{1}$ and $D_{2}$ in $2\mathrm{D}$. The corresponding narrow region
and its boundary parts in $2\mathrm{D}$ are defined by
\begin{align*}
\Omega_{r}&=\left\{(x_{1},x_{2})\in\mathbb{R}^{2}:~ -\frac{1}{2}h(x_{1})<x_{2}<\frac{1}{2}h(x_{1}),~ |x_{1}|<r\right\},\\
\Sigma^{+}_{r}&=\left\{(x_{1},x_{2})\in\mathbb{R}^{2}:~ x_{2}=\frac{1}{2}h(x_{1}),~ |x_{1}|<r\right\},\\
\Sigma^{-}_{r}&=\left\{(x_{1},x_{2})\in\mathbb{R}^{2}:~ x_{2}=-\frac{1}{2}h(x_{1}),~ |x_{1}|<r\right\},\\
\partial\Omega_{r}&=\left\{(x_{1},x_{2})\in\mathbb{R}^{2}:~ -\frac{1}{2}h(x_{1})<x_{2}<\frac{1}{2}h(x_{1}),~ |x_{1}|=r\right\},
\end{align*}
respectively. The mathematical representation for the problem in the narrow region $\Omega_{r}$ in $2\mathrm{D}$ keeps the same as $3\mathrm{D}$ except that ((\ref{11.10})(c)) should be substituted by
\begin{equation*}
\mathbf{u}=\mathbf{U}+\omega_{0}(\nu_{1}\mathbf{e}_{2}-\nu_{2}\mathbf{e}_{1}),\quad\quad\mathrm{on}\;\Sigma^{+}_{r}.
\end{equation*}

\subsection{Velocity split}

As for the case of $2\mathrm{D}$, $\boldsymbol{\nu}$ can be expressed as
\begin{equation*}
\boldsymbol{\nu}=x_{1}\mathbf{e}_{1}+\left(\frac{1}{2}|x_{1}|^{m}-R\right)\mathbf{e}_{2}.
\end{equation*}
Thus,
\begin{align}\label{9.01}
\mathbf{u}|_{\Sigma^{+}_{r}}=&\left[U_{1}-\omega_{0}\left(\frac{1}{2}|x_{1}|^{m}-R\right)\right]\mathbf{e}_{1}+(U_{2}+\omega_{0}x_{1})\,\mathbf{e}_{2}.
\end{align}
Using $(\ref{9.01})$ and $\mathbf{u}|_{\Sigma^{-}_{r}}=\mathbf{0}$, we rearrange the terms as follows:
\begin{align}
\mathbf{u}|_{\Sigma^{+}_{r}}=&\frac{1}{2}\left[\mathbf{U}+\omega_{0}\left(R-\frac{1}{2}|x_{1}|^{m},x_{1}\right)\right]+\frac{1}{2}(U_{1}+\omega_{0}R)\,\mathbf{e}_{1}+\frac{1}{2}U_{2}\mathbf{e}_{2}\notag\\
&-\frac{1}{4}\omega_{0}|x_{1}|^{m}\,\mathbf{e}_{1}+\frac{1}{2}\omega_{0}x_{1}\mathbf{e}_{2},\label{9.02}\\
\mathbf{u}|_{\Sigma^{-}_{r}}=&\frac{1}{2}\left[\mathbf{U}+\omega_{0}\left(R-\frac{1}{2}|x_{1}|^{m},x_{1}\right)\right]-\frac{1}{2}(U_{1}+\omega_{0}R)\,\mathbf{e}_{1}-\frac{1}{2}U_{2}\mathbf{e}_{2}\notag\\
&+\frac{1}{4}\omega_{0}|x_{1}|^{m}\,\mathbf{e}_{1}-\frac{1}{2}\omega_{0}x_{1}\mathbf{e}_{2}.\notag
\end{align}
By the linearity of the Stokes flow and (\ref{9.02}), we carry out a decomposition of the velocity $\mathbf{u}$ in $\Omega_{r}$ as follows:
\begin{equation}\label{111.15}
\mathbf{u}=
\sum\limits^{4}_{i=0}\mathbf{u}^{(i)},~~\quad\mathrm{in}\;2\mathrm{D},
\end{equation}
with $\mathbf{u}^{(i)},\,i\in\{0,\cdots,4\}$, solving problem $((\ref{11.10})\mathrm{(a)})$, $((\ref{11.10})\mathrm{(b)})$, $((\ref{11.10})(\mathrm{e}))$ and satisfying
\begin{flalign}\label{8.08}
&(0)~\mathbf{u}^{(0)}|_{\Sigma^{\pm}_{r}}=\frac{1}{2}\left[\mathbf{U}+\omega_{0}\left(R-\frac{1}{2}|x_{1}|^{m},x_{1}\right)\right],\notag\\
&(1)~\mathbf{u}^{(1)}|_{\Sigma^{\pm}_{r}}=\pm\frac{1}{2}(U_{1}+\omega_{0}R)\,\mathbf{e}_{1},\notag\\
&(2)~\mathbf{u}^{(2)}|_{\Sigma^{\pm}_{r}}=\pm\frac{1}{2}U_{2}\mathbf{e}_{2},\\
&(3)~\mathbf{u}^{(3)}|_{\Sigma^{\pm}_{r}}=\mp\frac{1}{4}\omega_{0}|x_{1}|^{m}\,\mathbf{e}_{1},\notag\\
&(4)~\mathbf{u}^{(4)}|_{\Sigma^{\pm}_{r}}=\pm\frac{1}{2}\omega_{0}x_{1}\mathbf{e}_{2}.\notag
\end{flalign}
Similarly as (\ref{11.17}), we get
\begin{align*}
\begin{cases}
\mathbf{F}=\sum\limits^{4}_{i=0}\int_{\Sigma^{+}_{r}}\boldsymbol{\sigma}(\mathbf{u}^{(i)})\mathbf{n}\,dS:=\sum\limits^{4}_{i=0}\mathbf{F}^{(i)},\\
\mathbf{T}=\sum\limits^{4}_{i=0}\int_{\Sigma^{+}_{r}}\boldsymbol{\nu}\times\boldsymbol{\sigma}(\mathbf{u}^{(i)})\mathbf{n}\,dS:=\sum\limits^{4}_{i=0}\mathbf{T}^{(i)}.
\end{cases}
\end{align*}
Similar explanations in the decomposition $(\ref{8.08})$ are omitted here. According to the case of $3\mathrm{D}$ in Section 3.3, similar assumptions for constructing the lubrication approximations for the fluid velocity $\mathbf{u}$ can be made for $2\mathrm{D}$.

\section{Proof of Theorem \ref{thm2}}\label{F4}
\begin{proof}[The proof of Theorem \ref{thm2}]
Its proof consists of the following six parts.

{\bf Part 1. Asymptotics of $\overline{\mathbf{F}}^{(1)}$ and $\overline{\mathrm{T}}^{(1)}$.}
Similarly as section 4.1, we let
\begin{align*}
\overline{\mathbf{u}}^{(1)}&=(U_{1}+\omega_{0}R)
\left(
\begin{array}{c}
          H(x_{1})x_{2} \\
          -A(x_{1})-B(x_{2})\frac{x_{2}^{2}}{2}
\end{array}
 \right),
\end{align*}
which together with the incompressibility and ((\ref{8.08})(1)) implies
\begin{equation*}
H=\frac{1}{h},\quad B=\partial_{1}H,\quad A=-B\frac{h^{2}}{8},
\end{equation*}
from which it follows by choosing $\overline{p}^{(1)}=C$ that
\begin{align}
\overline{\mathbf{F}}^{(1)}&=-\frac{2\mu(U_{1}+\omega_{0}R)\Gamma_{11}^{(m)}}{\varepsilon^{1-\frac{1}{m}}}\,\mathbf{e}_{1}+O(\mathbf{1}),\label{8.11}\\
\overline{\mathrm{T}}^{(1)}&=-\frac{2R\mu(U_{1}+\omega_{0}R)\Gamma_{11}^{(m)}}{\varepsilon^{1-\frac{1}{m}}}+O(1).\label{8.12}
\end{align}

{\bf Part 2. Asymptotics of $\overline{\mathbf{F}}^{(2)}$ and $\overline{\mathrm{T}}^{(2)}$.}
Compared with section 4.3, we set
\begin{align*}
\overline{\mathbf{u}}^{(2)}&=U_{2}
\left(
\begin{array}{c}
          -A_{1}(x_{1})-3B_{1}(x_{1})x_{2}^{2} \\
          A_{2}(x_{1})x_{2}+B_{2}(x_{1})x_{2}^{3}
\end{array}
 \right),
\end{align*}
which combining the incompressibility and $((\ref{8.08})(2))$ yields
\begin{equation*}
A_{1}=\frac{3}{2}\frac{x_{1}}{h},\quad B_{1}=-\frac{2x_{1}}{h^{3}},\quad A_{2}=\partial_{1}A_{1},\quad B_{2}=\partial_{1}B_{1}.
\end{equation*}
According to $(\ref{10.1})$, we choose
\begin{equation*}
\overline{p}^{(2)}=\mu U_{2}(-A_{2}+3B_{2}x_{2}^{2}-6G(x_{1})),
\end{equation*}
where
$$G(x_{1})=-\int^{x_{1}^{2}}_{r^{2}}\frac{1}{(\varepsilon+|t|^{\frac{m}{2}})^{3}}dt.$$

By a direct calculation, we deduce that
\begin{align}
\overline{\mathbf{F}}^{(2)}&=-\frac{12\mu U_{2}\Gamma_{33}^{(m)}}{\varepsilon^{3-\frac{3}{m}}}\,\mathbf{e}_{2}+O(\mathbf{1}),\quad\overline{\mathrm{T}}^{(2)}=0+O(1).\label{8.16}
\end{align}

{\bf Part 3. Asymptotics of $\overline{\mathbf{F}}^{(3)}$ and $\overline{\mathrm{T}}^{(3)}$.}
Choose
\begin{align*}
\overline{\mathbf{u}}^{(3)}&=\omega_{0}
\left(
\begin{array}{c}
          H(x_{1})x_{2} \\
          -A(x_{1})-B(x_{1})\frac{x_{2}^{2}}{2}
\end{array}
 \right).
\end{align*}
In light of the incompressibility and $((\ref{8.08})(3))$, we derive
\begin{equation*}
H=-\frac{1}{2}+\frac{\varepsilon}{2h},\quad B=\partial_{1}H,\quad A=-B\frac{h^{2}}{8}.
\end{equation*}

Let $\overline{p}^{(3)}=C$. Similarly, a direct calculation gives
\begin{align}
\overline{\mathbf{F}}^{(3)}&=\mathbf{0}+O(\mathbf{1}),\quad\overline{\mathrm{T}}^{(3)}=0+O(1).\label{8.19}
\end{align}

{\bf Part 4. Asymptotics of $\overline{\mathbf{F}}^{(4)}$ and $\overline{\mathrm{T}}^{(4)}$.}
Let
\begin{align*}
\overline{\mathbf{u}}^{(4)}&=\omega_{0}
\left(
\begin{array}{c}
          -A_{1}(x_{1})-3B_{1}(x_{1})x_{2}^{2} \\
          A_{2}(x_{1})x_{2}+B_{2}(x_{1})x_{2}^{3}
\end{array}
 \right).
\end{align*}
Due to the incompressibility and $((\ref{8.08})(4))$, we have
\begin{equation*}
A_{1}=\frac{3}{4}\frac{x^{2}_{1}}{h},\quad B_{1}=-\frac{x^{2}_{1}}{h^{3}},\quad A_{2}=\partial_{1}A_{1},\quad B_{2}=\partial_{1}B_{1}.
\end{equation*}

Similar to the construction given in section 4.6, we assume
\begin{equation*}
\overline{p}^{(4)}=\mu\omega_{0}(-A_{2}(x_{1})+3B_{2}(x_{1})x_{2}^{2}-6G(x_{1})),
\end{equation*}
where $\partial_{1}G=B_{1}$.

By definition, we have
\begin{align*}
\boldsymbol{\sigma}(\overline{\mathbf{u}}^{(4)})\mathbf{n}:=&\mu\omega_{0}(I_{1}\mathbf{e}_{1}+I_{2}\mathbf{e}_{2}),\\
\boldsymbol{\nu}\times\boldsymbol{\sigma}(\overline{\mathbf{u}}^{(4)})\mathbf{n}:=&\mu\omega_{0}J,
\end{align*}
where
\begin{align*}
I_{1}=&2(-\partial_{1}A_{1}-3(\partial_{1}B_{1})x_{2}^{2})n_{1}-\frac{\overline{p}^{(4)}}{\mu\omega_{0}}n_{1}+((\partial_{1}A_{2})x_{2}+(\partial_{1}B_{2})x_{2}^{3}-6B_{1}x_{2})n_{2},\\
I_{2}=&((\partial_{1}A_{2})x_{2}+(\partial_{1}B_{2})x_{2}^{3}-6B_{1}x_{2})n_{1}+2(A_{2}+3B_{2}x_{2}^{2})n_{2}-\frac{\overline{p}^{(4)}}{\mu\omega_{0}}n_{2},\\
J=&((\partial_{1}A_{2})x_{2}+(\partial_{1}B_{2})x_{2}^{3}-6B_{1}x_{2})n_{1}\nu_{1}+2(A_{2}+3B_{2}x_{2}^{2})\nu_{1}n_{2}-\frac{\overline{p}^{(4)}}{\mu\omega_{0}}\nu_{1}n_{2}+\frac{p^{(4)}}{\mu\omega_{0}}n_{1}\nu_{2}\\
&-2(-\partial_{1}A_{1}-3(\partial_{1}B_{1})x_{2}^{2})n_{1}\nu_{2}-((\partial_{1}A_{2})x_{2}+(\partial_{1}B_{2})x_{2}^{3}-6B_{1}x_{2})n_{2}\nu_{2}.
\end{align*}

Consider
\begin{align*}
G(x_{1})&=-\int^{x_{1}}_{-r}\frac{t^{2}}{(\varepsilon+|t|^{m})^{3}}\,dt.
\end{align*}

On one hand, it follows that
\begin{align*}
&\int_{\Sigma^{+}_{r}}-\frac{\overline{p}^{(4)}}{\mu\omega_{0}}n_{1}=-6(\varepsilon+r^{m})\Phi_{32}^{(m)}(r;\varepsilon)+6\Phi_{22}^{(m)}(r;\varepsilon)+O(1),\\
&\int_{\Sigma^{+}_{r}}((\partial_{1}A_{2})x_{2}+(\partial_{1}B_{2})x_{2}^{3}-6B_{1}x_{2})n_{2}=-6\Phi_{22}^{(m)}(r;\varepsilon)+O(1),\\
&\int_{\Sigma^{+}_{r}}-\frac{\overline{p}^{(4)}}{\mu\omega_{0}}n_{2}=12r\Phi_{32}^{(m)}(r;\varepsilon)+O(1),\\
&\int_{\Sigma^{+}_{r}}((\partial_{1}A_{2})x_{2}+(\partial_{1}B_{2})x_{2}^{3}-6B_{1}x_{2})n_{2}=-6\Phi_{22}^{(m)}(r;\varepsilon)+O(1).
\end{align*}
Thus, for $1<m\leq\frac{3}{2}$, we have
\begin{align}
\overline{\mathbf{F}}^{(4)}&=-\frac{3\mu\omega_{0}r\Gamma_{33}^{(m)}}{\varepsilon^{3-\frac{3}{m}}}\left(r^{m-1}\mathbf{e}_{1}-2\mathbf{e}_{2}\right)+O(\mathbf{1});\label{8.20}
\end{align}
for $m>\frac{3}{2}$,
\begin{align}
\overline{\mathbf{F}}^{(4)}=&-\frac{3\mu\omega_{0}\Gamma_{33}^{(m)}(r^{m}+\varepsilon)}{\varepsilon^{3-\frac{3}{m}}}\,\mathbf{e}_{1}+\frac{6\mu\omega_{0}r\Gamma_{33}^{(m)}}{\varepsilon^{3-\frac{3}{m}}}\,\mathbf{e}_{2}+O(\mathbf{1}).\label{8.24}
\end{align}

On the other hand, a direct calculation yields
\begin{align*}
&\int_{\Sigma^{+}_{r}}-((\partial_{1}A_{2})x_{2}+(\partial_{1}B_{2})x_{2}^{3}-6B_{1}x_{2})n_{2}\nu_{2}\\
=&-6R\Phi_{22}^{(m)}(r;\varepsilon)+\frac{9}{m}\Phi_{12}^{(m)}(r;\varepsilon)+O(1),\\
&\int_{\Sigma^{+}_{r}}((\partial_{1}A_{2})x_{2}+(\partial_{1}B_{2})x_{2}^{3}-6B_{1}x_{2})n_{1}\nu_{1}=9\Phi_{12}^{(m)}(r;\varepsilon)+O(1),\\
&\int_{\Sigma^{+}_{r}}\bigg(2(A_{2}+3B_{2}x_{2}^{2})\nu_{1}n_{2}-\frac{\overline{p}^{(4)}}{\mu\omega_{0}}\nu_{1}n_{2}\bigg)\\
=&\int^{r}_{-r}6x_{1}\int^{x_{1}}_{-r}\frac{t^{2}}{(\varepsilon+|t|^{m})^{3}}dtdx_{1}-6\Phi_{12}^{(m)}(r;\varepsilon)+O(1),
\end{align*}
and
\begin{align*}
&\int_{\Sigma^{+}_{r}}\frac{\overline{p}^{(4)}}{\mu\omega_{0}}n_{1}\nu_{2}\\
=&\frac{3}{2}\int^{r}_{-r}h\partial_{1}h\int^{x_{1}}_{-r}\frac{t^{2}}{(\varepsilon+|t|^{m})^{3}}dtdx_{1}-\frac{3\varepsilon+6R}{2}\int^{r}_{-r}\partial_{1}h\int^{x_{1}}_{-r}\frac{t^{2}}{(\varepsilon+|t|^{m})^{3}}dtdx_{1}.
\end{align*}
Therefore, we obtain that for $1<m\leq\frac{3}{2}$,
\begin{align}\label{3.001}
\overline{\mathrm{T}}^{(4)}&=\frac{3}{4}\frac{\mu\omega_{0}(4r^{2}+r^{2m}-4Rr^{m})\Gamma_{33}^{(m)}}{\varepsilon^{3-\frac{3}{m}}}+O(1);
\end{align}
for $\frac{3}{2}<m<\frac{5}{3}$,
\begin{align}\label{3.003}
\overline{\mathrm{T}}^{(4)}=&-\frac{3\mu\omega_{0}R\Gamma_{33}^{(m)}}{\varepsilon^{2-\frac{3}{m}}}+\frac{3}{4}\frac{\mu\omega_{0}(4r^{2}+r^{2m}-4Rr^{m})\Gamma_{33}^{(m)}}{\varepsilon^{3-\frac{3}{m}}}+O(1);
\end{align}
for $m=\frac{5}{3}$,
\begin{align}\label{3.004}
\overline{\mathrm{T}}^{(4)}=&-\frac{18}{5}\mu\omega_{0}|\ln\varepsilon|-\frac{3\mu\omega_{0}R\Gamma_{33}^{(\frac{5}{3})}}{\varepsilon^{\frac{1}{5}}}\notag\\
&+\frac{3}{4}\frac{\mu\omega_{0}(4r^{2}+r^{\frac{10}{3}}-4Rr^{\frac{5}{3}})\Gamma_{33}^{(\frac{5}{3})}}{\varepsilon^{\frac{6}{5}}}+O(1);
\end{align}
for $\frac{5}{3}<m<3$,
\begin{align}\label{3.005}
\overline{\mathrm{T}}^{(4)}=&-\frac{3\mu\omega_{0}R\Gamma_{33}^{(m)}}{\varepsilon^{2-\frac{3}{m}}}-\frac{3\mu\omega_{0}\Gamma_{35}^{(m)}}{\varepsilon^{3-\frac{5}{m}}}\notag\\
&+\frac{3}{4}\frac{\mu\omega_{0}(4r^{2}+r^{2m}-4Rr^{m})\Gamma_{33}^{(m)}}{\varepsilon^{3-\frac{3}{m}}}+O(1);
\end{align}
for $m=3$,
\begin{align}\label{3.006}
\overline{\mathrm{T}}^{(4)}=&\frac{3}{2}\mu\omega_{0}|\ln\varepsilon|-\frac{3\mu\omega_{0}R\Gamma_{33}^{(3)}}{\varepsilon}-\frac{3\mu\omega_{0}\Gamma_{35}^{(3)}}{\varepsilon^{\frac{4}{3}}}\notag\\
&+\frac{3}{4}\frac{\mu\omega_{0}(4r^{2}+r^{6}-4Rr^{3})\Gamma_{33}^{(3)}}{\varepsilon^{2}}+O(1);
\end{align}
for $m>3$,
\begin{align}\label{3.007}
\overline{\mathrm{T}}^{(4)}=&\frac{\mu\omega_{0}}{2m}\frac{(18+3m)\Gamma_{13}^{(m)}+6m\Gamma_{23}^{(m)}}{\varepsilon^{1-\frac{3}{m}}}-\frac{3\mu\omega_{0}R\Gamma_{33}^{(m)}}{\varepsilon^{2-\frac{3}{m}}}\notag\\
&-\frac{3\mu\omega_{0}\Gamma_{35}^{(m)}}{\varepsilon^{3-\frac{5}{m}}}+\frac{3}{4}\frac{\mu\omega_{0}(4r^{2}+r^{2m}-4Rr^{m})\Gamma_{33}^{(m)}}{\varepsilon^{3-\frac{3}{m}}}+O(1).
\end{align}

{\bf Part 5. Asymptotics of $\overline{\mathbf{F}}^{(0)}$ and $\overline{\mathrm{T}}^{(0)}$.}
Suppose
\begin{align*}
\overline{\mathbf{u}}^{(0)}&=\frac{1}{2}
\left(
\begin{array}{c}
          U_{1}+\omega_{0}(R-\frac{|x_{1}|^{m}}{2}) \\
          U_{2}+\omega_{0}x_{1}
\end{array}
 \right).
\end{align*}
Choose $\overline{p}^{(0)}=C$. It is easy to verify that
\begin{align}
\overline{\mathbf{F}}^{(0)}&=\mathbf{0}+O(\mathbf{1}),\quad\overline{\mathrm{T}}^{(0)}=0+O(1).\label{8.27}
\end{align}

{\bf Part 6. Justification.}
Similarly as in section 4.8, combining with (\ref{8.11})--(\ref{8.27}), we complete the proof of Theorem \ref{thm2}.

\end{proof}

As for the particles with partially flat boundary, we also have similar results in $2\mathrm{D}$. A direct application of Theorem \ref{thm3} gives
\begin{theorem}\label{thm5}
Assume that $D_{1},D_{2}\subset\Omega\subseteq\mathbb{R}^{2}$ are defined as above, condition (\ref{2.001}) holds. Let $\mathbf{u}\in \mathbf{H}^{1}(\Omega)$ be the solution of problem \eqref{ZAD123}. Then, for a sufficiently small $\varepsilon$,
\begin{align*}
F_{1}=&-\left(\alpha_{\frac{1}{2}}\varepsilon^{-\frac{1}{2}}+\alpha_{1}\varepsilon^{-1}+\alpha_{\frac{3}{2}}\varepsilon^{-\frac{3}{2}}+\alpha_{2}\varepsilon^{-2}+\alpha_{\frac{5}{2}}\varepsilon^{-\frac{5}{2}}+\alpha_{3}\varepsilon^{-3}\right)+O(1),\\
F_{2}=&-\left(\beta_{\frac{3}{2}}\varepsilon^{-\frac{3}{2}}+\beta_{2}\varepsilon^{-2}-\beta_{\frac{5}{2}}\varepsilon^{-\frac{5}{2}}+\beta_{3}\varepsilon^{-3}\right)+O(1),\\
\mathrm{T}=&-\left(\gamma_{0}\ln\varepsilon+\gamma_{\frac{1}{2}}\varepsilon^{-\frac{1}{2}}+\gamma_{1}\varepsilon^{-1}+\gamma_{\frac{3}{2}}\varepsilon^{-\frac{3}{2}}+\gamma_{2}\varepsilon^{-2}+\gamma_{\frac{5}{2}}\varepsilon^{-\frac{5}{2}}+\gamma_{3}\varepsilon^{-3}\right)+O(1),
\end{align*}
where
\begin{align*}
\alpha_{\frac{1}{2}}=&2\mu(U_{1}+\omega_{0}R)\Gamma_{11}^{(2)}+3\mu\omega_{0}\Gamma_{33}^{(2)},~\;\alpha_{1}=2\mu(U_{1}+\omega_{0}R)s+3\mu\omega_{0}s,\\
\alpha_{\frac{3}{2}}=&\mu\omega_{0}\big(3(r-s)^{2}\Gamma_{33}^{(2)}+3s^{2}\Gamma_{31}^{(2)}\big),~\; \alpha_{2}=\mu\omega_{0}\big(3s(r-s)^{2}+2s^{3}\big),\\ \alpha_{\frac{5}{2}}=&3\mu\omega_{0}s^{2}(r-s)^{2}\Gamma_{31}^{(2)},~\;\alpha_{3}=2\mu\omega_{0}s^{3}(r-s)^{2},~\;\beta_{\frac{3}{2}}=6\mu(2U_{2}-\omega_{0}r)\Gamma_{33}^{(2)},\\
\beta_{2}=&6\mu s(U_{2}-\omega_{0}r),~\;\beta_{\frac{5}{2}}=6\mu\omega_{0}s^{2}r\Gamma_{31}^{(2)},~\;\beta_{3}=4\mu s(3U_{2}r^{2}-U_{2}s^{2}-\omega_{0}rs^{2}),\\
\gamma_{0}=&\frac{9}{2}\mu\omega_{0}s,~\;\gamma_{1}=\mu\omega_{0}s(3R-s^{2}+6)+2\mu R(U_{1}+\omega_{0}R)s,\\
\gamma_{\frac{1}{2}}=&2\mu R(U_{1}+\omega_{0}R)\Gamma_{11}^{(2)}+\frac{3}{4}\mu\omega_{0}\big(-6s^{2}\Gamma_{11}^{(2)}-2s^{2}\Gamma_{21}^{(2)}+s^{2}\Gamma_{31}^{(2)}+4R\Gamma_{33}^{(2)}+4\Gamma_{35}^{(2)}\big),\\
\gamma_{\frac{3}{2}}=&\frac{1}{4}\mu\omega_{0}\big[\big(12R(r-s)^{2}-3(r-s)^{4}-12r^{2}\big)\Gamma_{33}^{(2)}+12Rs^{2}\Gamma_{31}^{(2)}+36s^{2}\big],\\
\gamma_{2}=&\frac{1}{4}\mu\omega_{0}\big(8s^{3}R-12sr^{2}+24s^{3}+12Rs(r-s)^{2}-3s(r-s)^{4}\big),\\
\gamma_{\frac{5}{2}}=&\frac{1}{4}\mu\omega_{0}(-12s^{2}r^{2}+12s^{4}+12Rs^{2}(r-s)^{2}-3s^{2}(r-s)^{4})\Gamma_{31}^{(2)},\\
\gamma_{3}=&\frac{1}{2}\mu\omega_{0}\big(-20r^{2}s^{3}+12s^{5}+20Rs^{3}(r-s)^{2}-5s^{3}(r-s)^{4}\big),
\end{align*}
and the remainder terms of order $O(1)$ depend only on $\mu$, $R$, $r$, $s$, $U_{i}$ and $\omega_{0}$, but not on $\varepsilon$.
\end{theorem}
The proof of Theorem \ref{thm5} is similar to Theorem \ref{thm2} and thus omitted here.

\noindent{\bf{\large Acknowledgements.}} Li was partially supported by NSFC (11631002, 11971061) and BJNSF (1202013).

\bibliographystyle{plain}

\def\cprime{$'$}

\end{document}